\documentclass[12pt]{article}
\usepackage{fancyhdr,amsmath, graphicx, psfrag, color ,amsfonts, layout}

\font\fontauthors=cmcsc10 scaled \magstep1

 at 8truept
 at 8truept

\newtheorem{Th}{Theorem}
\newtheorem{Def}[Th]{Definition}
\newtheorem{Cor}[Th]{Corollary}
\newtheorem{Lem}[Th]{Lemma}
\newtheorem{Prop}[Th]{Proposition}
\newtheorem{Rem}[Th]{Remark}

\newcommand{\pff}{\noindent {\sc Proof.}\ }

\def\QED{\hfill$\square$ \vskip 10pt}

\newtheorem{prerem}[Th]{Remark}

\def\convepsilon{\smash{\mathop{\longrightarrow}\limits _{\varepsilon\to 0}}}

\def\convcont{\smash{\mathop{\longrightarrow}\limits _{t\to \infty}}}

\def\egalLoi{{~\mathop{= }\limits^{\rond L}}~}

\def\Sp{\mathop{\rm Sp}\nolimits}
\def\Span{\mathop{\rm Span}\nolimits}
\def\Supp{\mathop{\rm Supp}\nolimits}
\def\stirling2 #1#2{\left\{\begin{matrix} #1\\#2\end{matrix}\right\}}

\def\Var{\mathop{\rm Var}\nolimits}

\newcommand{\R}{\rm I\!R}
\newcommand{\C}{\;{}^{{}_\vert}\!\!\!{\rm C}}
\newcommand{\N}{{{\rm I}\!{\rm N}}}

\newcommand{\Q}{{\rm Q}\kern-.65em{}^{{}_/}}

\def\1{{\bf 1}}

\def\g#1{\mathbb #1}
\def\rond#1{\mathcal #1}
\def\R{\g R}
\def\C{\g C}
\def\N{\g N}
\begin{document}

%~
%\vskip 15pt
\begin{center}
\LARGE{\bf {Limit distributions for\\

multitype branching processes of\\

{\LARGE{$m$-ary} }search trees \footnote{
{\it 2000 Mathematics Subject Classification.} Primary: 60C05. Secondary: 60J80, 05D40.

{\it Key words and phrases.}
Martingale.
Characteristic function.
Embedding in continuous time.
Multitype branching process.
Smoothing transformation. Absolute continuity. Support. Exponential moments.
}}}\\
\end{center}

\begin{center}
by
\end{center}

\begin{center}
{\fontauthors
Brigitte Chauvin,
%\footnote{INRIA Rocquencourt, project Algorithms and Laboratoire de Math\'ematiques de Versailles, CNRS, UMR 8100 --  INRIA Domaine de Voluceau B.P.105, 78153 Le Chesnay CEDEX (France).},
Quansheng Liu and
Nicolas Pouyanne
%\footnote{Laboratoire de Math\'ematiques de Versailles, CNRS, UMR 8100 -- Universit\'e de Versailles - St-Quentin,  45 avenue des Etats-Unis, 78035 Versailles CEDEX (France).
}
\end{center}

\begin{center}
{\it 1 December 2011}\\
\end{center}
\noindent{\bf Abstract.} A particular continuous-time multitype branching process is considered, it is the continuous-time embedding of a discrete-time process which is very popular in theoretical computer science: the $m$-ary search tree ($m$ is an integer). There is a well-known phase transition: when $m\leq 26$, the asymptotic behavior of the process is Gaussian, but for $m\geq 27$ it is no more Gaussian and a limit $W$ of a complex-valued martingale arises. Thanks to the branching property it appears as a solution of a {\it smoothing} equation of the type
$Z \egalLoi e^{ -\lambda T  }  (Z^{(1)} + \cdots + Z^{(m)})$,
where  $\lambda\in \g C$,  the $Z^{(k)}$ are independent copies   of $Z$ and $T$ is a $\g R_+$-valued random variable, independent of the $Z^{(k)}$.
This distributional equation is extensively studied by various approaches. The existence and unicity of  solution of the equation are proved by contraction methods. The fact that the distribution of  $W$ is absolutely continuous  and  that its support is the whole complex plane is shown via Fourier analysis.
%%
%%
%% Q: it does not sound very well to say "the determination ... which is ...;
%% also: "its Fourier transform " does not seem very clear. Aussi, 3 fois "existence" ne sonne pas tres bien.
%%
%%Analysis of its Fourier transform provides the existence of a density for $W$ together with
% the determination of the support of $W$ which is the whole complex plane.
%%
%%
Finally, the existence of exponential moments of $W$ is obtained by considering $W$ as the limit of a complex Mandelbrot cascade.

\bigskip

\tableofcontents

%%%%%%%%%%%%%%
%\section{Introduction}
%%%%%%%%%%%%%%

%\setcounter{equation}{0}
\setcounter{section}{0}

\section{Introduction}\label{intro}

\renewcommand{\theequation}{\thesection.\arabic{equation}}
\setcounter{equation}{0}

Consider a continuous-time multitype branching process $(X(t), t\geq 0)$. Types are seen as colors of particles and there are $m-1$ colors, where $m\geq2$ is an integer. The reproduction of the process is given by a particular matrix $R$ (written in (\ref{matrix})), and any particle of colour
$j$ lives a random time of exponential distribution with parameter $j$. Such a classical process is considered for  example in
Athreya and Ney \cite{AN} or Janson \cite{Jan} and it is precisely defined in Section 2.

When it is stopped  at the $n$-th jump time, this process is nothing but the composition vector process $(X^{DT}_n, n\geq 0)$ say, of an $m$-ary search tree, which is an important algorithmic structure in computer science.
A numerous literature is devoted to the asymptotic behavior of this composition vector.
A famous phase transition appears. When $m\leq 26$, the random vector admits a central limit theorem
with convergence in distribution to a Gaussian vector: see Mahmoud and Pittel \cite{MahPit}, Lew and Mahmoud \cite{LewMah}.
When  $m\geq 27$, it has been proved that
\begin{equation*}%\label{discrete}
X_n^{DT} = nv_1 +  \Re(n^{\lambda_2} W^{DT} v_2 ) + o(n^{\sigma_2}) \quad a.s.,
\end{equation*}
where $\lambda_2$ is a complex number having a real part in $] \frac 12, 1[$ (it is an eigenvalue of the replacement matrix $R$),
where $v_1, v_2$ are deterministic vectors and $W^{DT}$ is the limit of a complex-valued martingale.
Heated conjectures  about the  random variable $W^{DT}$ remain open
% appearing in the expansion of this vector
(see \cite{ChaPou04}, \cite{Pou05}, Chern and Hwang \cite{ChernHwang}, Mahmoud \cite{Mahmoud},  Janson \cite{Jan}).
%it is complex valued, $o(\ )$ holds almost surely  and in any $L^p$-space ($p\geq 1$), and the moments of $W^{DT}$ can be recursively calculated.

This article is focused on the asymptotic behavior of the continuous-time process $(X(t))$.
The links between continuous-time and discrete-time processes are detailed in Section 4.
In particular, we give the so-called \emph{martingale connection} that relates the almost sure limits of both processes.
%This work thus constitutes a step towards a better knowledge of $W^{DT}$.

\bigskip
Inspired by the methods used for a two-color P\'olya urn in \cite{ChaPouSah}, we first prove in Section 3 that for $m\geq 27$, $X(t) $ admits the following  asymptotic expansion:
$$
X(t) =  e^t \xi v_1 (1+o(1)) + \Re(e^{\lambda_2 t} W v_2) (1+o(1)) \; \mbox{ a.s. and in } L^p \;\; \forall p\geq 1,
$$
where  $\xi$ is a Gamma distributed random variable and $W$ a $\C$-valued one.

We are interested in the limit random variables  $W_k, k= 1, \dots , m-1$, each corresponding to $X_k(t)$ which denotes the process $X(t)$  when it starts from one particle of color $k$. Using the branching property, a system of dislocation equations is written for the  random vectors $X_k(t)$ in Section \ref{dislocationX}.
A system of fixed point equations satisfied by the corresponding limit laws is then derived in Section \ref{fixedW}.
In particular, the complex-valued random variable $W_1$ is a solution of the fixed point equation
\begin{equation}\label{fixed}
Z \egalLoi e^{ -\lambda_2 T  }  (Z^{(1)} + \cdots + Z^{(m)}),
\end{equation}
where  $\{Z^{(k)}: k\geq 1 \} $ are independent copies   of $Z$,   $T = \tau_{(1)} + \dots + \tau_{(m-1)}$,  $\{ \tau_{(j)}: j\geq 1\}$ are random variables independent of each other and  independent of $\{Z^{(k)}\}$,
each  $\tau_{(j)} $ has distribution  ${\cal E}xp(j)$ (we denote by ${\cal E}xp(j)$ the exponential distribution of parameter $j$: it has density $x\mapsto je^{-jx}$ on $]0,\infty[ $).

\medskip
Further properties of $W_1$ are derived from a fine study of Eq.~(\ref{fixed}).
We first show in Theorems~\ref{W} and \ref{contractiond2*} that Eq.~(\ref{fixed}) admits a
unique square-integrable solution having a given mean.
In particular, this implies that Eq.~(\ref{fixed}) characterizes the distribution of $W_1$.
This result is proven by two contraction methods applied to the corresponding smoothing transformations.
The first one deals with suitable spaces of probability measures where the classical Wasserstein metric is
adapted to the complex field;
it leads to Theorem~\ref{W}.
The second contraction method, that gives a proof for Theorem~\ref{contractiond2*}, consists in working on
Fourier transforms of solutions and provides a somehow simpler proof.
Furthermore, this second method gives a result of existence and unicity for solutions of the convolution equation
\begin{equation*}
\Phi(t) = \int_0^{+\infty} \Phi^m(t-u) f_T(u) du, \qquad t\in \C,
\end{equation*}
in a convenient space of functions, where $f_T$ denotes the density of $T$
(see Remark~\ref{remEqConvolution}).

\medskip
Once the characterization of $W_1$ by Eq.~(\ref{fixed}) is proven, it suffices to derive properties of
solutions of this distributional equation.
We show in this way the following results on the law of $W_1$.

\begin{Th}
\label{s1th1}
When $m\geq 27$, the complex-valued random variable $W_1$ admits a density and its support is the whole
complex plane.
Its Fourier transform satisfies
$$
\g Ee^{i\langle t,W_1\rangle}=O(|t|^{-a})
$$
when $|t|\to +\infty$, for some $a>1$.
\end{Th}

This theorem is a direct consequence of Theorem~\ref{density} that provides such properties for solutions
of~(\ref{fixed}) admitting a nonzero mean.
Our proof consists in showing successively that the characteristic function of any solution has modulus equal
to~$1$ only at the origin, that it tends to zero at infinity, and finally that it is of order $O(|t|^{-a})$ as $|t|\rightarrow \infty$ for some
$a>1$ so that it is square-integrable on $\g C$.
In the approach we need to prove a  non-lattice property of Eq.~(\ref{fixed}) \emph{via} Gelfand-Schneider
theorem, using the algebraicity of $\lambda _2$
(see proof of Lemma~\ref{lemma0}).

\begin{Th}
\label{s1th2}
When $m\geq 27$, the random variable $W_1$ admits exponential moments in a neighbourhood of the origin
of the complex plane.
If $L_1(z)=\g Ee^{zW_1}$ denotes its Laplace series, then $L_1$ is holomorphic near $0$ and, after a change
of variable, the function $z\mapsto -\frac\rho zL_1\left( z^{-\lambda _2}\right)$ is a solution of the differential
equation
$$y^{(m-1)}=y^m.
$$
\end{Th}

Theorem \ref{s1th2} is  immediately derived from Theorems~\ref{Expmoments} and
\ref{thSerieLaplace} just as Theorem \ref{s1th1} was derived from Theorem~\ref{density}.
To prove Theorems \ref{Expmoments} and
\ref{thSerieLaplace}, we consider a solution of~(\ref{fixed}) as the limit of a complex Mandelbrot cascade.
The results are a consequence of fine analytical properties of the Fourier transform of the limit variable.

\medskip
The paper is organized as follows.

The continuous-time multitype branching process is defined in Section~\ref{model}.
Its relation with the $m$-ary search tree
%space requirements composition vector %
is detailed in
Section~\ref{connection},  while Section~\ref{martingaleConnexion} is devoted to the second order asymptotic
expansion of $(X(t))_{t\geq 0}$ and to its connection with the corresponding
discrete process.
In Section~\ref{secDislocation}, we use the branching property of the process to show that the martingale limits
of the continuous-time process are related by a system of equations in law so that the fixed point
equation~(\ref{fixed}) emerges.
These first four sections constitute the first part of the paper.

The second part of the paper consists in putting the focus on Eq.~(\ref{fixed}) that turns out to characterize
the distribution of $W_1$ so that all results on solutions provide results on $W_1$.
In Section~\ref{smoothing} we define the natural smoothing transform associated with Eq.~(\ref{fixed})
and we show that it defines a contraction in the space of square-integrable probability measures with given
mean.
Results on the support and on absolute continuity of solutions are obtained in Section~\ref{secDensite}.
Finally, Section~\ref{laplace} is devoted to the exponential moments and the Laplace series of solutions.

%%%%%%%%%%%%%%%%%%%%
\section{Definition of the branching process} %$X = (X(t), t\geq0)$}
\label{model}
%%%%%%%%%%%%%%%%%%%%%%%
%%%%
%% B,N : on trouvait pas tres beau d'avoir des maths dans un titre de section
%%%%%%
\renewcommand{\theequation}{\thesection.\arabic{equation}}
\setcounter{equation}{0}

In this section we introduce the definition of the continuous time multitype branching process $ (X(t))$, and present the spectral decomposition of its transition matrix.

\subsection{Infinitesimal generator}
\label{generator}

In the whole paper, the underlying vector space is ${\g R}^{m-1}$ or sometimes
${\g C}^{m-1}$.
Let $R$ be the following square matrix of order $m-1$:
\begin{equation}\label{matrix}
R=
\left(
\begin{array}{ccccc}
-1&1&&&\\
&-1&1&&\\
&&-1&\ddots&\\
%&&\ddots&\ddots&\\
&&&\ddots&1\\
m&&&&-1
\end{array}
\right)
,
\end{equation}
and for $k= 1, \dots , m-1$, let $w_k$  be the $k$-th  row vector of $R$:
when $1\leq k\leq m-2$, the $k$-th coordinate of $w_k$ equals $-1$,  the $(k+1)$-th equals $1$ and all the others are $0$;
 $w_{m-1}$ has $m$ as first coordinate, $-1$ as last one, and $0$ for all others.

Let $G$ be the operator defined on functions $f$ from ${\g C}^{m-1}$ to any real or complex vector space
by the following formula:
for any vector $v$ in ${\g C}^{m-1}$,
\begin{equation}
\label{defPhi}
G(f)(v) = \sum_{k=1}^{m-1} k l_k(v) [f(v+w_k) - f(v) ],
\end{equation}
where  $l_k$ are the coordinate forms:  $l_k(x_1,\dots ,x_{m-1})=x_k$.

\begin{Def}
The right-continuous process $X = (X(t), t\geq0)$ is the only continuous time Markov process with state
space $\g R^{m-1}$ having $G$ as infinitesimal generator.
\end{Def}

\medskip

\noindent
Equivalently, $X$ is a continuous time multitype branching process with $m-1$ types (or colors),
having $R$ as reproduction matrix.
The $k$-th coordinate of the vector $X(t)$, namely $l_k(X(t))$, is  the number of  particles of color $k$ at time $t$.
A particle of color $k$ lives a random exponential time with parameter $k$;  when it dies, it reproduces one particle of color $k+1$ if $k= 1, \dots , m-2$, and
 $m$ particles of color $1$ if $k=m-1$.

\medskip

This branching continuous time process can be thought as the embedded process of a discrete Markov chain
$X^{DT} = (X_n^{DT})_{n\in\g N}$ which is a P\'olya-type discrete Markov chain associated with
the node process of an $m$-ary search tree, an important algorithmic structure in computer science.
This connection is detailed in Section~\ref{connection}.

%%%%%%%%%%%%%%%%%%%%%
\subsection{Spectral decomposition}
%%%%%%%%%%%%%%%%%%%%%

Let $R_G$ be the matrix of $G$'s restriction to linear forms in the canonical basis $(l_k)_{1\leq k\leq m-1}$.
One immediately checks that
\begin{equation*}%\label{matrixG}
R_G=
\left(
\begin{array}{ccccc}
-1&1&&&\\
&-2&2&&\\
&&-3&\ddots&\\
%&&\ddots&\ddots&\\
&&&\ddots&m-2\\
m(m-1)&&&&-(m-1)
\end{array}
\right),
\end{equation*}
where an empty entry means a zero entry.
It has been established in many papers -- see for example Mahmoud \cite{Mahmoud}, Chern and Hwang
\cite{ChernHwang} or \cite{ChaPou04} -- and it can be easily checked that
$R_G$'s (unitary) characteristic polynomial is
\begin{equation}
\label{polycar}
\chi _{R_G}(\lambda )=
\prod _{k=1}^{m-1}(\lambda  +k)-m!
=\frac{\Gamma (\lambda +m)}{\Gamma (\lambda +1)}-m! ,
\end{equation}
where $\Gamma$ denotes Euler's Gamma function.
All eigenvalues are simple, $1$ being the one having the largest real part.

\begin{center}
\fbox{
\begin{minipage}{140 truemm}
\it
In the whole paper, $\lambda _2$ will denote $\chi _{R_G}$'s root having the second largest real part and a positive
imaginary part.
\end{minipage}
}
\end{center}
The famous phase transition on $m$-ary search trees already mentioned in the introduction is due to the fact
that
$$ \Re (\lambda _2)>1/2 \; \mbox{ if and only if } \;   m\geq 27. $$
See for example  \cite{ChaPou04}.
The assumption $\Re (\lambda _2)>1/2$  will be frequently used in the sequel.

%\noindent
We adopt the following notations:
\begin{equation}\label{notations}
\left\{
\begin{array}{l}
\displaystyle
\forall n\in\g Z_{\geq 0},~\binom{z}{n}=\frac{\Gamma (z+1)}{n!\Gamma (z-n+1)}=\frac{z(z-1)\dots (z-n+1)}{n!};\\
\displaystyle
%\forall z\notin \g Z_{\leq -1},~
H_m(z)=\sum _{1\leq k\leq m-1}\frac 1{z+k};\\
\displaystyle
u_1 (x_1,\dots ,x_{m-1})=\sum _{1\leq k\leq m-1}kx_k ;\\
\displaystyle
u_2 (x_1,\dots ,x_{m-1})=
%\sum _{1\leq k\leq m-1}\frac{\Gamma (\lambda _2+k)}{\Gamma (\lambda _2+1)\Gamma (k)}x_k;\\
\sum _{1\leq k\leq m-1}\binom{\lambda _2+k-1}{k-1}x_k;\\
\displaystyle
v_1=\frac 1{H_m(1)}\left( \frac 1{k(k+1)}\right) _{1\leq k\leq m-1};\\
%=\frac 1{H_m(1)}^t\left( \frac 1{1\times 2},\frac 1{2\times 3},\dots ,\frac1{(m-1)m}\right);\\
\displaystyle
v_2=\frac 1{H_m(\lambda _2)}
\left( \frac 1{k\binom {\lambda _2+k}{k}}\right) _{1\leq k\leq m-1}.
%=\frac 1{H_m(\lambda _2)}^t\left( \frac 1{\lambda _2+1},\dots ,\frac 1{(m-1)\binom {\lambda _2+m-1}{m-1}}\right).
\end{array}
\right.
\end{equation}
The linear forms $u_1$ and $u_2$ are eigenvectors of $G$, namely
$G (u_1)=u_1$ and $G (u_2)=\lambda _2u_2$.
The vectors $v_1$ and $v_2$ are left eigenvectors of $R_G$, respectively associated with the eigenvalues $1$ and $\lambda_2$.
They satisfy $u_1(v_1)=u_2(v_2)=1$ and $u_1(v_2)=u_2(v_1)=0$.
These eigendata had already been essentially computed in \cite{ChaPou04} and \cite{Pou05}.
Note that, since $\lambda _2$ is not real, $u_2$ and $v_2$ have nonreal coordinates.
%%%%%%%%%%%%%%%%%%%%%%%%%%%%%%%%%%%%%%%%%%%%%%%%%%%%%%%%%%%%%%%%%%%%%%%%%%%%%%%%%%%%%%%%%
%%%%
%%%% the notations $\overline{u_2}$ and $\overline{v_2}$ will be used for complex conjugacy.
%%%% Q: this need to be mentioned by the end of the intro.
%%%%%
%%%%%%%%%%%%%%%%%%%%%%%%%%%%%%%%%%%%%%%%%%%%%%%%%%%%%%

\section{ $m$-ary search trees and embedding}
\label{connection}

\renewcommand{\theequation}{\thesection.\arabic{equation}}
\setcounter{equation}{0}

%%%% B, N : on veut quand meme citer Janson, pour ne pas pretendre qu'on est les premiers ˆ ecrire cela

In this section we present the connection between $m$-ary search trees and multitype branching processes. This example of embedding of a discrete time process into a continuous time process has already been evoked by Janson in \cite{Jan}.

\subsection{$m$-ary search trees}
\label{mary_def}

We define here a discrete time Markov chain $X^{DT} = (X_n^{DT}, n \geq 0)$ with values in $\g N^{m-1}\setminus \{ 0\}$. The $i$-th coordinate of $X_n^{DT}$ is denoted by $X_n^{(i)}$ and has a ``physical'' meaning detailed hereafter. The Markov chain $X^{DT}$ is a random walk defined by an initial
vector $X^{DT}_0$ in $\g N^{m-1}\setminus \{ 0\}$  and by the following transition probabilities:

\noindent $\forall v\in \g N^{m-1}\setminus \{ 0\}$, $\forall k= 1, \dots , m-1$,
\begin{equation}
\label{transitionproba}
q(v, v+w_k) = \frac{kl_k(v)}{\sum_{j=1}^{m-1} jl_j(v)},
\end{equation}
where the increment vectors $w_k$ are given in Section  \ref{generator} and $l_k(v)$ denotes the $k$-th coordinate of the  vector $v$.

\bigskip

Classically (see Norris \cite{Norris} and for a synthetic exposition Bertoin \cite{Bertoin}), this discrete time  Markov chain is embedded in continuous time using a ``Poissonization'' of the time: given $X^{DT} $, one can recover $X=(X(t), t\geq 0)$ as follows. At time $0$, $X(0) =  X^{DT}_0$. For any vector $v\in\g R^{m-1}$, define\footnote{Note that $q=u_1$ where $u_1$ was defined by (\ref{notations}).}
$$
q(v):= \sum_{k=1}^{m-1} k l_k(v).
$$
Let $\tau_1$ be a random time exponentially distributed with parameter $q(X^{DT}_0)$. For any time $t\in [0, \tau_1[$, let $X(t) = X(0) = X^{DT}_0$. At time $\tau_1$, $X$ jumps from $v= X(0)$ to $v+w_k$ with probability given by formula (\ref{transitionproba}). More generally, let $\tau_0 = 0$ and for any $n\geq 1$, define the $n$-th jumping time $\tau_n$ by
$$
\tau_n = \sum_{i=0}^{n-1} \frac{\epsilon_i}{q(X_i^{DT})},
$$
where $\epsilon_i$ are independent random variables having the same exponential distribution with parameter $1$.
Let
$$
X(t) = X(\tau_n) = X_n^{DT}, \hskip 5mm \forall t\in [\tau_n, \tau_{n+1}[.
$$
At time $\tau_{n+1}$, $X$  jumps from $v=X(\tau_n)$ to $v+w_k$ with probability given by formula (\ref{transitionproba}). It is easy  to see that this embedded process $X(t)$ is the same one as the branching process defined in Section  \ref{generator}.

\bigskip

When $X^{DT}_0 = (1, 0, \dots , 0)$, each $X_n^{(i)}, i = 1, \dots , m-1$, can be seen as the number of nodes of type $i$ in a tree $T_n$: the sequence $(T_n, n\geq 0)$ is a sequence of random $m$-ary trees which grow  by successive insertions of keys in their leaves. Each node of these trees contains at most $m-1$ keys. Keys are i.i.d. random variables $x_i, i\geq 1$, with any diffusive distribution on the interval $[0,1]$. The tree $T_n, n\geq 0$, is recursively defined as follows: $T_0$ is reduced to an empty node-root; $T_1$ is reduced to a node-root which contains $x_1$, $T_2$ is reduced to a node-root which contains $x_1$ and $x_2$, ... , $T_{m-1}$ has a node-root containing $x_1, \dots x_{m-1}$.   As soon as the ($m-1$)-th key is inserted in the root, $m$ empty subtrees of the root are created, corresponding from left to right to the $m$ ordered intervals $I_1 = ]0,x_{(1)}[, \dots, I_m=]x_{(m-1)}, 1[$, where $0< x_{(1)} < \dots < x_{(m-1)} <1$ are the ordered $m-1$ first keys. Each following key $x_m, \dots $ is recursively inserted in the subtree corresponding to the unique interval $I_j$ to which it belongs. As soon as a node is saturated, $m$ empty subtrees of this node are created.

For each $i=\{1,\dots ,m-1\}$ and $n\geq 1$, $X_n^{(i)}$ is the number
of nodes in $T_n$ which contain $i-1$ keys (and $i$ gaps or free places) after insertion of the $n$-th key; such nodes are named nodes of type $i$. We don't worry about  the number of saturated nodes. The vector $X_n^{DT}$ is called  the composition vector of the $m$-ary search tree. It provides a model for the space requirement of the algorithm. One can refer to Mahmoud's book \cite{Mahmoud} for further details on search trees.

Notice that, in this dynamics, the insertion of a new key is {\it uniform} on the gaps; it can be read on the transition probabilities (\ref{transitionproba}).

\subsection{Embedding}

The embedding properties  are summarized in the following lemma.

\begin{Lem}\label{embed}
\

1) For any $n\geq 1$, the distribution of $\tau_n - \tau_{n-1}$ is ${\cal E}xp(n-1+N_0)$, where $N_0$ is the number of free places in $X(0)$: $N_0 = u_1(X(0))$.

\medskip

2) the processes $(\tau_n)_{n\geq 1}$ and $(X(\tau_n))_{n\geq 1}$ are independent.

\medskip

3) the processes $(X(\tau_n))_{n\geq 1}$ and $(X_n^{DT})_{n\geq 1}$ have the same distribution.
 \end{Lem}
\pff
Part 1) is a consequence of the fact that the minimum of $k$ independent ${\cal E}xp(1)$-distributed random variables is ${\cal E}xp(k)$-distributed,  and that the total number of free places at time $\tau_n$ equals $n-1+N_0$.

Part 2) is the classical independence between the jump chain and the jump times in such Markov processes. The initial states and evolution rules of both Markov chains in discrete time and in continuous time are the same ones, so that  Part 3) holds.~\QED

\noindent\underline{Convention}.  From now on, thanks to Part 3) of Lemma \ref{embed}, we will as usual suppose that the  discrete-time process and the continuous-time process are built on the \emph{same}
probability space on which
\begin{equation}
\label{embedding}
(X(\tau_n))_{n\geq 1} = (X_n^{DT})_{n\geq 1} \ \ a.s..
\end{equation}

\medskip

\noindent\underline{Remark}. The important benefit we get with the embedding is the independence in the
continuous-time process. This independence is the key point for the dislocation equations later on.

%%%%%%%%%%%%%%%%%%%%%%%%%%%%
\section{Asymptotics and martingale connection}
\label{martingaleConnexion}
%%%%%%%%%%%%%%%%%%%%%%%%%%%%%

In this section we present an order 2 expansion of the continuous time multitype branching process $(X(t))$ and its connection with the discrete time process
$(X_n^{DT})$ defined in Section \ref{mary_def}.

\renewcommand{\theequation}{\thesection.\arabic{equation}}
\setcounter{equation}{0}

\subsection{Asymptotics of the continuous time branching process}

With the notations of Section~\ref{model} and especially the formulae~(\ref{notations}), the random vector $X(t)$
admits the following order 2 expansion as $t$ goes to infinity.

\begin{Th}
\label{thAsymptotiqueCT}
{\bf (Asymptotics of continuous time process)}

\noindent
Suppose that $m\geq 27$.
Then, as $t$ tends to infinity,
\begin{equation}
\label{asymptotiqueCT}
X(t)=e^{t}\xi v_1\big( 1+\varepsilon_1(t)\big) +2\Re\left(e^{\lambda_2t}W v_2\right)\big( 1+\varepsilon_2(t)\big) +
%e^{\lambda_2t}
{\boldsymbol\varepsilon}_3 (t) ,
\end{equation}
where

\noindent
$\bullet$
$\xi$ is a positive {\rm Gamma}-distributed random variable with expectation $N_0= u_1(X(0))$
(total weighted number of particles at time $0$),

\noindent
$\bullet$
$W$ is a complex-valued random variable that admits moments of all orders $p\geq 1$ and
whose expectation equals $u_2(X(0))$,

\noindent
$\bullet$
the real-valued random variables $\varepsilon_1 (t)$ and $\varepsilon_2 (t)$ tend to $0$  as $t$ tends to
$+\infty$, almost surely and in any ${\rm L}^p$-space, $p\geq 1$,

\noindent
$\bullet$
the random vector $\boldsymbol\varepsilon_3 (t)$ is $o\left( e^{\lambda_2t}\right)$
%$\varepsilon_2 (t)$ is vector-valued and tends to $0$
as $t$ tends to $+\infty$, almost surely and in any ${\rm L}^p$-space, $p\geq 1$.
\end{Th}

\begin{center}
\fbox{
\begin{minipage}{140 truemm}
\it
In the whole paper, $W$ denotes our hero, namely the limit complex-valued random variable of the second
order term in $X(t)$'s expansion, as in Theorem~\ref{thAsymptotiqueCT}.
\end{minipage}
}
\end{center}

\begin{Rem}
One can reformulate (\ref{asymptotiqueCT}) as follows:

\noindent
in any basis of the form $(v_1,\Re(v_2), \Im(v_2),\dots )$,

\noindent
$\bullet$
$X(t)$'s first coordinate has the expansion $e^{t}\xi + o(e^t)$,

\noindent
$\bullet$
$X(t)$'s component on $\Span _{\g R}\left(\Re(v_2),  \Im(v_2)\right)$ has the expansion
$2\Re\left(e^{\lambda_2t}Wv_2 \right) + o(e^{\lambda_2t}) $,

\noindent
$\bullet$
all other coordinates are $o(e^{\lambda_2t})$.

\end{Rem}
\noindent {\sc Proof of Theorem}~\ref{thAsymptotiqueCT}.
Denote $\rond A$ the endomorphism of $\g R^{m-1}$ having $\, ^t\!R_G$ as matrix in the canonical basis. Let also $M(t)=\exp (-t\rond A)X(t)$, for any $t\geq 0$.
By standard arguments from multitype branching process theory,
$(M(t))_{t\geq 0}$ is a vector-valued martingale.
Since $m\geq 27$, the real part of $\lambda _2$ belongs to $]1/2,1[$ so that the projected martingales
$u_1(M(t))$ and $u_2(M(t))$ converge in ${\rm L}^p$ for any $p\geq 1$.
For proofs of these results, see for example Athreya and Ney \cite{AN} or Janson \cite{Jan}
(especially Lemma 10.2 of Janson's paper for the ${\rm L}^p$-boundedness, $X$ being here an
\emph{irreducible} process in the sense of~\cite{Jan}).
The random variables $\xi$ and $W$ are respectively defined by
\begin{equation}
\label{defxiW}
\left\{
\begin{array}{l}
\displaystyle
\xi =\lim _{t\to +\infty}u_1\left( e^{-t\rond A}X(t)\right) =\lim _{t\to +\infty}e^{-t}u_1\left( X(t)\right),\\ \\
\displaystyle
W=\lim _{t\to +\infty}u_2\left( e^{-t\rond A}X(t)\right) =\lim _{t\to +\infty}e^{-\lambda _2t}u_2\left( X(t)\right) .
\end{array}
\right.
\end{equation}
An alternative proof of the ${\rm L}^p$ convergence can be made using the techniques of \cite{Pou08},    as
developed in~\cite{ChaPouSah} for two-colour urn processes.
In particular, $\xi$'s distribution is attained by explicit computation of its moments:
for any nonnegative integer $p$, an elementary computation shows directly from~(\ref{defPhi}) that the
(so-called \emph{reduced}) polynomial
$$
Q:=u_1\left( u_1+1\right)\left( u_1+2\right)\dots \left( u_1+p-1\right)
$$
is an eigenvector for $X$'s infinitesimal generator $G$, associated with the eigenvalue~$p$.
Thus $\g EQ(X(t))=e^{pt} Q(X(0))$ for any $t$.
Besides, because of~(\ref{defxiW}), $Q(X(t))=e^{pt}\xi ^p(1+o(1))$ as $t$ tends to infinity, almost surely and
in ${\rm L^1}$.
Finally, the last two equalities provide
$$
\g E\xi^p=Q(X(0))=\frac{\Gamma \left( N_0+p\right)}{\Gamma \left( N_0\right)}.
$$
This shows that  the law of  $\xi$ is a Gamma distribution with parameter $N_0$
 since a Gamma distribution is completely  determined by its moments.
The matrix  $R_G$ is diagonalizable on $\g C$ since all roots of its characteristic polynomial are simple
(see~(\ref{polycar})).
Extending notations~(\ref{notations}), let $(u_\lambda)_{\lambda\in\Sp (\rond A)}$ be a basis of linear forms,
each $u_\lambda$ being an eigenform of $G$ associated with the (complex) eigenvalue $\lambda$.
Let also $(v_\lambda)_{\lambda\in\Sp (\rond A)}$ be the dual basis of $(u_\lambda)_{\lambda\in\Sp (\rond A)}$,
each $v_\lambda$ being thus a vector that satisfies $u_\lambda (v_\mu )=\delta _{\lambda,\mu }$
(Kronecker's notation).
Note that one can choose $u_{\lambda _2}=u_2$ and, consequently, $v_{\lambda _2}=v_2$
(cf. notations~(\ref{notations})).

\noindent
For any $t\geq 0$, split the spectral decomposition of the vector $X(t)$ with respect to $G$ into
four terms:
$$
X(t)=\sum _{\lambda\in\Sp (\rond A)}u_\lambda (X(t)).v_\lambda
=X_1(t)+X_2(t)+X_3(t)+X_4(t),
$$
where
\begin{equation*}
\left\{
\begin{array}{l}
X_1(t)=u_1(X(t))v_1, \\ \\
X_2(t)=u_{\lambda _2}(X(t))v_{\lambda _2}+u_{\overline{\lambda _2}}(X(t))v_{\overline{\lambda _2}}, \\ \\
\displaystyle
X_3(t)=\sum _{1/2<\Re\lambda <\Re\lambda _2}u_\lambda (X(t))v_\lambda, \\ \\
X_4(t)=\sum _{\Re\lambda <1/2}u_\lambda (X(t))v_\lambda.
\end{array}
\right.
\end{equation*}
Note that this partition of $\Sp (\rond A)$ is valid because $\frac 12$ is not an eigenvalue of $\rond A$ as can be
checked from~(\ref{polycar}).
We deal separately with these four components of $X(t)$.
Define $\boldsymbol\varepsilon _3$ by $\boldsymbol\varepsilon _3(t)=X_3(t)+X_4(t)$, for any $t\geq 0$.

\vskip 5pt
$\bullet$
The formulae~(\ref{defxiW}) provide directly the asymptotics
$$
\left\{
\begin{array}{l}
X_1(t)=(e^t\xi +o(e^t))v_1, \\ \\
X_2(t)=(e^{\lambda _2t}W+o(e^{\lambda _2t}))v_2+\overline{(e^{\lambda _2t}W+o(e^{\lambda _2t}))v_2}
=2\Re\left((e^{\lambda_2t}W +o(e^{\lambda_2t}))v_2\right) ,
\end{array}
\right.
$$
leading to the first two terms of the expansion~(\ref{asymptotiqueCT}).

\vskip 5pt
$\bullet$
Suppose that $\lambda$ is an eigenvalue of $\rond A$ such that $\frac 12<\Re\lambda <\Re\lambda _2$.
Then, with the same general arguments as in the very beginning of the proof, it can be seen that
$$
u_\lambda (M(t))=e^{-t\lambda}u_\lambda (X(t))
$$
and that $(u_\lambda (M(t)))_{t\geq 0}$ is a convergent
martingale, bounded in any ${\rm L}^p$, $p\geq 1$.
In particular, $u_\lambda (X(t)) =  o(e^{\lambda _2t})$ as $t$ tends to infinity, almost surely and in any
${\rm L}^p$, $p\geq 1$.
This shows that $X_3(t)$ is $o(e^{\lambda _2t})$ when $t\to +\infty$.

\vskip 5pt
$\bullet$
It remains to deal with the small eigenvalues, namely with all $\lambda$ such that $\Re\lambda <\frac 12$.
\begin{Lem}
\label{lemmePetitesVp}
Suppose that $\lambda$ is an eigenvalue such that $\Re\lambda <\frac 12$ and let $\eta >0$.
Then, $e^{-(\frac 12 +\eta)t}u_{\lambda}(X(t))$ is bounded almost surely and in any
${\rm L}^p$-space, $p\geq 1$.
\end{Lem}
The proof of this lemma is given just hereafter.
Therefore, if $\Re\lambda <\frac 12$, then
$$
e^{-\lambda _2t}u_\lambda (X(t))
=e^{(1/2+\eta -\lambda _2)t}\left[ e^{-(\frac 12 +\eta)t}u_{\lambda}(X(t))\right] \convcont \ 0
$$
almost surely as soon as $0<\eta <\Re\lambda _2-\frac 12$.
Such $\eta$ exist because $\Re\lambda _2>\frac 12$.
This shows that $X_4(t)$ is $o(e^{\lambda _2t})$ when $t\to +\infty$.
The same argument holds for the ${\rm L}^p$ convergence, making the proof complete.
\QED

\noindent {\sc Proof of Lemma~\ref{lemmePetitesVp}}.
The main idea consists in taking advantage of the following fact:
when $t$ belongs to the interval $[\tau_n, \tau_{n+1}[$, the vector $X(t)$ remains equal to $X_n^{DT}$.
This being considered, we make use of the moment bounds of the discrete time process that can be found in
\cite{Pou08} (Theorem 3.4 (1)):
when $\Re\lambda <\frac 12$,
\begin{equation}
\label{momentdiscret}
\forall p\geq 1,\forall \varepsilon >0,
\hskip 10pt  \g E |u_{\lambda}(X_n^{DT})|^p = O\left( n^{p(\frac 12 + \varepsilon)}\right) ,
\hskip 30pt n\rightarrow + \infty.
\end{equation}

$\bullet$ Almost sure bound:
we prove that
\begin{equation}
\label{asbound}
\lim_{C\rightarrow +\infty} \g P\left(\exists t>0,  \ e^{-(\frac 12 +\eta)t}\left| u_{\lambda}(X(t))\right| >C\right) = 0,
\end{equation}
which suffices to get the almost sure boundedness.
Let $C>0$, $\eta >0$ and let $\lambda$ be an eigenvalue such that $\Re\lambda <\frac 12$.
The jump time $\tau _n$ tends almost surely to $+\infty$ which is a classical result that can be deduced
from Lemma~\ref{embed}, so that
$$
\g P\left( \exists t>0,  \ e^{-(\frac 12 +\eta)t} |u_{\lambda}(X(t))| >C\right)
\leq\displaystyle
\sum_{n\geq 0} \g P\left( \exists t\in [\tau_n, \tau_{n+1}[ ,  \ e^{-(\frac 12 +\eta)t} |u_{\lambda}(X(t))| >C\right).
$$
Since $X(t)=X_n^{DT}$ for any $t\in [\tau_n, \tau_{n+1}[$, this leads to
$$
\g P\left( \exists t>0,  \ e^{-(\frac 12 +\eta)t} |u_{\lambda}(X(t))| >C\right)
\leq\displaystyle
\sum_{n\geq 0}\g P\left(\left| u_{\lambda}\left( X_n^{DT}\right)\right| >Ce^{(\frac 12 +\eta)\tau_n}\right)
$$
Conditioning with respect to $\tau _n$, using Markov inequality and the fact that $\tau_n$ and $X_n^{DT}$ are
independent, one gets successively, for any $p\geq 1$:
$$
\begin{array}{rl}
\g P\left( \exists t>0,  \ e^{-(\frac 12 +\eta)t} |u_{\lambda}(X(t))| >C\right)&
\leq\displaystyle
\sum_{n\geq 0}
\g E\left( \g P\left(\left| u_{\lambda}(X_n^{DT})\right| > C e^{(\frac 12 +\eta)\tau_n}\big| \tau_n\right) \right)
\\ \\
&\leq\displaystyle
\sum_{n\geq 0}
\g E\left( \frac{ \g E\left| u_{\lambda}\left( X_n^{DT}\right)\right| ^p}{C^p e^{p(\frac 12 +\eta)\tau_n}} \right)
\\ \\
&=\displaystyle
\frac 1{C^p}\sum_{n\geq 0}
\g E\left| u_{\lambda}\left( X_n^{DT}\right)\right| ^p\g E\left( e^{-p(\frac 12 +\eta)\tau_n} \right).
\end{array}
$$
The density of the $n$-th jump time $\tau_n$ is the function
$$
u\in\g R\longmapsto ne^{-u}\left( 1- e^{-u}\right) ^{n-1} \1_{\R_+}(u),
$$
so that its Laplace transform can be elementarily computed: for any $s\geq 0$,
$$
\g E (e^{-s\tau_n})
=\frac{n! \Gamma(s+1)}{\Gamma(s+1+n)} \sim \Gamma(s+1) n^{-s}, \hskip 1cm n\rightarrow +\infty.
$$
Together with (\ref{momentdiscret}), this leads to:
$\forall \eta >0, \forall \varepsilon >0$, $\forall p\geq 1$,
$$
\g E\left| u_{\lambda}\left( X_n^{DT}\right)\right| ^p\g E\left( e^{-p(\frac 12 +\eta)\tau_n} \right)
=O\left(\frac{1}{n^{p(\eta - \varepsilon)}}\right), \hskip 1cm n\rightarrow +\infty
$$
which is the general term of a convergent series as soon as one takes $\varepsilon <\eta$ and
$p>\frac 1{\eta -\varepsilon}$.
Finally, letting $C$ tend to infinity shows~(\ref{asbound}).

\vskip 5pt
$\bullet$ Bound in ${\rm L}^p$-space:
let $p\geq 1$ and $t>0$.
Then,
$$
\left\lVert e^{-(\frac 12+\eta )t}u_\lambda \left( X(t)\right)\right\rVert _p^p
=e^{-(\frac 12+\eta )pt}\g E\left| u_\lambda \left( X(t)\right)\right| ^p.
$$
Using the relation with the discrete time process $(X_n^{DT})_n$, one has successively
$$
\begin{array}{rl}
\left\lVert e^{-(\frac 12+\eta )t}u_\lambda \left( X(t)\right)\right\rVert _p^p&
=\displaystyle e^{-(\frac 12+\eta )pt}\sum _{n\geq 0}
\g E\Big( \1 _{\tau _n\leq t<\tau _{n+1}}\big| u_\lambda \left( X(t)\right)\big| ^p\Big)\\
&=\displaystyle e^{-(\frac 12+\eta )pt}\sum _{n\geq 0}
\g E\Big( \1 _{\tau _n\leq t<\tau _{n+1}}\left| u_\lambda \left( X_n^{DT}\right)\right| ^p\Big)\\
&=\displaystyle e^{-(\frac 12+\eta )pt}
\sum _{n\geq 0}\g E\Big( \1 _{\tau _n\leq t<\tau _{n+1}}\Big)
\g E\Big(\left| u_\lambda \left( X_n^{DT}\right)\right| ^p\Big) ,
\end{array}
$$
where the last equality holds due to the independence between $\tau _n$ and $X_n^{DT}$. Besides, $\tau _n$ and $\tau _{n+1}-\tau _n$ are independent and $\tau _{n+1}-\tau _n$ is $\rond Exp (n+N_0)$-distributed (see~(\ref{embed})), so that, using the density of $\tau _n$ written above, one gets
$$
\begin{array}{rl}
\g E\Big( \1 _{\tau _n\leq t<\tau _{n+1}}\Big)&
=\g E\Big(\1 _{t\geq \tau _n}\g E\Big( \1 _{\tau _{n+1}-\tau _n\geq t-\tau _n}\big|\tau _n\Big)\Big)\\
&=\g E\Big(\1 _{t\geq \tau _n}e^{-(n+N_0)(t-\tau _n)}\Big)\\
&\displaystyle=\int _0^te^{-(n+N_0)(t-u)}ne^{-u}\left( 1-e^{-u}\right)^{n-1}du\\
&\displaystyle\leq ne^{-(n+1)t}\int _0^t\left( e^u-1\right) ^{n-1}e^{u}du
=\left( 1-e^{-t}\right)^ne^{-t}.
\end{array}
$$
Thus,
$$
\left\lVert e^{-(\frac 12+\eta )t}u_\lambda \left( X(t)\right)\right\rVert _p^p
\leq e^{-t}e^{-(\frac 12+\eta )pt}
\sum _{n\geq 0}\left( 1-e^{-t}\right)^n
\g E\Big(\left| u_\lambda \left( X_n^{DT}\right)\right| ^p\Big).
$$
Let now $\varepsilon >0$.
On one hand, (\ref{momentdiscret}) implies that
$$
\g E\Big(\left| u_\lambda \left( X_n^{DT}\right)\right| ^p\Big)
= O\left( n^{p(\frac 12 + \varepsilon)}\right) .
$$
On the other hand, Stirling's formula applied to generalized binomial coefficients yields classically
that for any $\alpha\in\g C$,
$$
[z^n](1-z)^{-\alpha -1}
%=(-1)^n\cnp{-\alpha -1}{n}
=\frac{n^{\alpha}}{\Gamma (\alpha +1)}\left( 1+O\left(\frac 1n\right)\right) ,
$$
where the notation $[z^n] A(z)$ means the  coefficient of $z^n$ in the power expansion of $A(z)$ at the origin. Consequently,
$$
\g E\Big(\left| u_\lambda \left( X_n^{DT}\right)\right| ^p\Big)
=O\left( [z^n](1-z)^{-1-p(\frac 12+\varepsilon)}\right) .
$$
This implies that for any $\varepsilon >0$, there exists a constant $C_\varepsilon$ such that for any $t>0$,
$$
\left\lVert e^{-(\frac 12+\eta )t}u_\lambda \left( X(t)\right)\right\rVert _p^p
\leq C_\varepsilon e^{-t}e^{-(\frac 12+\eta )pt}
\left( 1-\left( 1-e^{-t}\right)\right)^{-1-p(\frac 12+\varepsilon)}
=C_\varepsilon e^{-pt(\eta - \varepsilon )}.
$$
It suffices to take $\varepsilon =\eta/2$ to conclude that the ${\rm L}^p$-norm of
$e^{-(\frac 12+\eta )t}u_\lambda \left( X(t)\right)$ is bounded above.
\QED

\begin{Rem}
The distribution of $W$ is infinitely divisible, because it is the limit of infinitely divisible ones, obtained by scaling and projection of infinitely divisible ones. Indeed, in finite time, for any $x_0 \in \g R^{m-1}$, denote by $(X_{x_0}(t), t\geq 0)$ the process $(X(t),  t\geq 0)$ defined in Section \ref{generator} starting from initial state $x_0$. By the branching property
$$
X_{x_0}(t) \egalLoi [n] X_{\frac{x_0}n}(t),
$$
where the notation $[n]X$ denotes the sum of $n$ independent copies of the random variable $X$. This fact has already been noticed by Janson (\cite{Jan}, proof of Theorem~3.9).
\end{Rem}

\subsection{Martingale connection}

In this subsection, we use the embedding equality (\ref{embedding}) to deduce connections between the asymptotic behaviours of $X_n^{DT}$ when $n\rightarrow + \infty$ and $X(t)$ when $t\rightarrow + \infty$.

\medskip

For $m\geq 27$, it has been proved in \cite{ChaPou04} and \cite{Pou05} that
\medskip

\begin{equation}\label{discrete}
X_n^{DT} = nv_1 +  \Re(n^{\lambda_2} W^{DT} v_2 ) + o(n^{\sigma_2}) \;\; \mbox{ a.s. and in } L^p,  \;\; \forall p\geq 1,
\end{equation}
where $v_1, v_2$ are the deterministic vectors defined in (\ref{notations}),
$W^{DT}$ is a complex-valued martingale limit,
% the $o(\ )$ is almost sure  and in any $L^p$-space, $p\geq 1$
and
the moments of $W^{DT}$ can be recursively calculated.

\medskip
\begin{Prop}
\label{propMC} The following two assertions hold:
\begin{equation}
\label{PropEQ1}
\lim_{ n\rightarrow +\infty} n e^{ -\tau_n} = \xi \; \mbox{ a.s. and in }  \; L^p, \; \forall p\geq 1,
 \end{equation}
\begin{equation}
\label{martingaleconnection}
W = \xi^{\lambda_2}\  W^{DT} \;\; \mbox{ a.s. with }\;  \xi \mbox{ and }  W^{DT} \mbox{ independent}.
\end{equation}
\end{Prop}
The equality (\ref{martingaleconnection}) is commonly referred to as ``martingale connection''.

\medskip

\pff
We first prove (\ref{PropEQ1}).
 Applying the first projection to the embedding equality  (\ref{embedding}), we obtain that
$$
u_1( X(\tau_n)) = u_1( X_n^{DT}))\ \ a.s.,
$$
where $u_1$ has been defined in (\ref{notations}). This is the total number of free places at time $\tau_n$, and is  equal to
 $n-1+N_0 = n(1+o(1))$.  Therefore, by (\ref{defxiW}) and the fact that the splitting times
$\tau_n$ tend almost surely to $+\infty$ when $n$ goes to $+\infty$, we have
$$
\xi= \lim_{t\rightarrow +\infty} e^{-t}u_1 ( X(t) ) = \lim_{ n\rightarrow +\infty} n e^{ -\tau_n}\ \ a.s..
$$
This gives  (\ref{PropEQ1}).

We then prove (\ref{martingaleconnection}).  Applying the second projection to the embedding equality  (\ref{embedding}) we obtain
$$
u_2( X(\tau_n)) = u_2( X_n^{DT}))\ \ a.s.,
$$
where $u_2$ has been defined in (\ref{notations}). Using again  (\ref{defxiW}) and the fact that $\tau_n$ goes to $+\infty$ when $n$ goes to $+\infty$,
we get
$$
W= \lim_{t\rightarrow +\infty} e^{-\lambda_2 t}u_2 ( X(t) ) = \lim_{ n\rightarrow +\infty}  e^{ -\lambda_2 \tau_n}u_2( X_n^{DT}) .
$$
Therefore  (\ref{martingaleconnection})  follows from (\ref{PropEQ1})
 and the asymptotics in discrete time given in (\ref{discrete}). \QED
\begin{Rem}
Fill and Kapur (\cite{FillKapur}) proved that $W^{DT}$ is the unique solution in the space of probability distributions with a given mean and finite second absolute moment of the fixed point equation
\begin{equation}\label{DTfixedpoint}
Z \egalLoi \sum_{k=1}^m (V_k)^{\lambda_2} Z^{(k)},
\end{equation}
where the $V_k$ are the spaces in the statistical order of $(m-1)$ i.i.d. random variables uniformly distributed  on $[0,1]$.
Because $\R$-valued stable distributions are solutions of the fixed point equation (\ref{DTfixedpoint}) when $\lambda$ is a real number, it is somehow natural to ask whether a $\lambda$-stable distribution is a solution of Eq. (\ref{DTfixedpoint}) for a complex number $\lambda$. By $\lambda$-stable we mean operator-stable when the operator is given by a two dimensional matrix $\lambda = \sigma + i\tau =
%\label{matrix}
\left(
\begin{array}{cc}
\sigma&-\tau\\
\tau&\sigma
\end{array}
\right)
$ as introduced by Sharpe \cite{Sharpe}. It is known since Hudson et al. \cite{HVW} that a $\lambda$-stable distribution has infinite moments of order $p$ for $p > 1/\Re(\lambda)$. Consequently, neither $W$ nor $W^{DT}$ (which have moments of any order) can be stable distributions.
\end{Rem}

%%%%%%%%%%%%%%%%%%%
\section{A distributional equation}
\label{secDislocation}
%%%%%%%%%%%%%%%%%%%

\renewcommand{\theequation}{\thesection.\arabic{equation}}
\setcounter{equation}{0}

In this section we derive a distributional equation satisfied by  the limit variable of the continuous-time branching process with an appropriate norming.
We shall see that this equation characterizes the limit distribution.

\subsection{Vectorial finite time dislocation equations}\label{dislocationX}

Let us write dislocation equations for the continuous-time branching process at finite time $t$.
We write  $X_j(t)$ for $X(t)$ when the process starts from $X(0) = e_j$, where $e_j$ denotes the  $j$-th vector of the canonical basis of $\g R^{m-1}$ (whose $j$-th component is $1$ and all the others are $0$). This means that the process starts from an ancestor of type $j$.

Notice that the first splitting time $\tau_1$ changes of distribution depending on the ancestor's type; denote by $\tau_{(j)}, j=1, \dots , m-1$,  the first splitting time when the process starts from $X(0) = e_j$. Thus $\tau_{(j)}$ is ${\cal E}xp(j)$ distributed.

The branching property applied at the first splitting time gives:

\begin{equation}
\label{dislocationfinitetime}
\forall t>\tau_1,~
\left\{
\begin{array}{l}
\displaystyle
X_1(t) \egalLoi X_2(t-\tau_{(1)}), \\ \\
X_2(t) \egalLoi X_3(t-\tau_{(2)}), \\ \\
\dots \\ \\
X_{m-2}(t) \egalLoi X_{m-1}(t-\tau_{(m-2)}), \\ \\
X_{m-1}(t) \egalLoi [m] X_1(t-\tau_{(m-1)}),
\end{array}
\right.
\end{equation}
where the notation $[m]X$ denotes the sum of $m$ independent copies of the random variable $X$.

\medskip

% \noindent\underline{Notation}:
In the following, we denote
\begin{equation}\label{notationT}
T = \tau_{(1)} + \dots + \tau_{(m-1)},
\end{equation}
where the $\tau_{(j)}$ are independent of each other and each
 $\tau_{(j)}$ is ${\cal E}xp(j)$  distributed.
Let us give some elementary properties of $T$ that we shall need.
It is easy to see that $T$ has density
\begin{equation}
\label{densitefT}
f_{T} (u) = (m-1) e^{-u} (1- e^{-u})^{m-2}\1_{\R_+}(u), \quad u\in \R,
\end{equation}
so that $e^{-T}$ has a Beta distribution with parameters  $1$ and $m-1$.
A straightforward change of variable under the integral
shows that for any complex number $\lambda$ such that
$\Re (\lambda )>0$,
\begin{eqnarray}
\label{beta}
  \g E e^{-\lambda T} & = & \int _0^{+\infty}e^{-\lambda u}f_T(u)du =(m-1)B(1+\lambda ,m-1)   \\
                   & = & \frac{ (m-1)! } { \prod_{k=1}^{m-1} (\lambda + k)},
\end{eqnarray}
where $B$ denotes Euler's Beta function:
\begin{equation*}
 B(x,y) = \int _0^1 u^{x-1}(1-u)^{y-1} du = \frac{ \Gamma (x) \Gamma (y)}{\Gamma (x+y)},
 \quad \Re x >0, \Re y >0.
 \end{equation*}
In particular,
\begin{equation}
\label{T}
 m  \g E | e^{-\lambda T}| \;
   \left\{
\begin{array}{l}
\displaystyle
 <1 \mbox { if }  \Re (\lambda) >1,   \\
 = 1 \mbox { if }  \Re (\lambda) = 1,   \\
 > 1 \mbox { if }  \Re (\lambda) <1.  \\
\end{array}
\right.
\end{equation}

\subsection{Distributional equation satisfied by the limit variable}
\label{fixedW}
%\subsection{Fixed point equations for $W$}\label{fixedW}

After projections of variables $X_j(t)$ (the process starts  from $X(0) = e_j$) with $u_2$, scaling with $e^{ -\lambda_2 t}$ and taking the limit when $t$ goes to infinity, we get the variables
$$
W_j:= \lim_{t\rightarrow +\infty} e^{-\lambda_2 t}u_2 ( X_j(t) ),
$$
so that the system (\ref{dislocationfinitetime}) on $X_j(t)$ leads to the following system of distributional equations on   $W_j$:
\begin{equation}
\label{system}
\left\{
\begin{array}{l}
\displaystyle
W_1\egalLoi e^{ -\lambda_2\tau_{(1)}}W_2 , \\ \\
W_2\egalLoi e^{ -\lambda_2\tau_{(2)}}W_3 , \\ \\
\dots \\ \\
W_{m-2}\egalLoi e^{ -\lambda_2\tau_{(m-2)}}W_{m-1} , \\ \\
W_{m-1}\egalLoi e^{ -\lambda_2\tau_{(m-1)}}[m] W_1 .
\end{array}
\right.
\end{equation}
This shows that $W_1$ is a solution of  the following fixed point equation:
\begin{equation}
\label{fixedpointmary}
Z \egalLoi e^{ -\lambda_2 T } (Z^{(1)} + \cdots + Z^{(m)}),
\end{equation}
where $ Z^{(i)} $ are independent copies of $Z$, which are also independent of $T$.

In terms of  the Fourier transform
$$ \varphi(t) :=  \g E \exp \{ i \langle t,Z \rangle \}   = \g E \exp \{ i \Re (\overline t Z) \},    \quad t \in \C,    $$
where  $ \langle x,y\rangle   =   \Re (\overline x y ) = \Re (x) \Re (y)  + \Im (x) \Im (y) $,
the equation (\ref{fixedpointmary}) reads
\begin{equation}
\label{fixedpointmaryF}
\varphi (t)
      =  \int_0^{+\infty} \varphi^m (t e^{-\overline{\lambda_2} u}) f_T(u) du, \quad t\in \C,
      \end {equation}
where  $f_T $ is defined by (\ref{densitefT}).
Notice that this functional equation can also be written in a convolution form: if
$\Phi(t):=\varphi (e^{\overline{\lambda_2}t})$   for any $t\in \C$, then
  $\Phi$ satisfies the following functional equation:
\begin{equation}
\label{psi}
\Phi(t) = \int_0^{+\infty} \Phi^m(t-u) f_T(u)du, \quad t\in \C.
\end{equation}

In the following sections, we prove that the distributional equation (\ref{fixedpointmary}) characterizes the law of $W_1$ and we get several results on $W_1$: for example  $W_1$ has a density on the whole complex plane, and admits exponential moments. All these results appear as a particular case of a slightly more general situation given hereafter. From now on, for any complex number $\lambda$, consider the distributional equation
\begin{equation}
\label{fixedpoint}
Z \egalLoi e^{ -\lambda T } (Z^{(1)} + \cdots + Z^{(m)}),
\end{equation}
where $ Z^{(i)} $ are independent copies of $Z$, which are also independent of $T$. On Fourier transforms, it reads
\begin{equation}
\label{fixedpointF}
\varphi (t)
      =  \int_0^{+\infty} \varphi^m (t e^{-\overline{\lambda} u}) f_T(u) du, \quad t\in \C,
      \end {equation}
where  $f_T $ is defined by (\ref{densitefT}).

Notice that when Z is a solution of the distributional equation (\ref{fixedpoint}), with finite and non zero first moment, then $\lambda $ is a root of the polynomial function~(\ref{polycar}). In particular, $\lambda$ is an algebraic number.

%%%%%%%
%\section{The smoothing transformation}
\section{The smoothing transformation}
\label{smoothing}

\renewcommand{\theequation}{\thesection.\arabic{equation}}
\setcounter{equation}{0}

A solution of the distributional equation (\ref{fixedpoint}) is a fixed point of the associated smoothing
transformation defined hereafter by~(\ref{defK}).
Endowing a suitable space of probability measures with two distances, we prove that the smoothing
transformation is a contraction for both metrics.
This provides two alternative approaches for the study of Eq.~(\ref{fixedpoint}) by contraction method.
Using the Wasserstein distance as a first metric, we adapt the classical contraction method developed
in~\cite{Gui90}, \cite{Ros92} and~\cite{RR}.
The second metric is defined in terms of Fourier transforms of measures;
it provides a short proof of our result.

\vskip 5pt
For any complex number $C$, let  ${\cal M}_2(C)$ be the space of probability distributions on $\C$ admitting a second absolute moment and having $C$ as expectation.

Let $\lambda$ be a complex number.
For any probability measure $\mu$ on $\g C$, let
\begin{equation}
\label{defK}
K \mu :=   \rond L \big( e^{ -\lambda T} (Z^{(1)} + \cdots +  Z^{(m)})  \big),
\end{equation}
where $T$ is given by (\ref{notationT}), $Z^{(i)} $ are independent random variables of law $\mu$,  which are
also independent of $T$.
Following Durrett and Liggett \cite{DL83} who considered the case of real random variables,  we call $K$  the
\emph{smoothing transformation}.
Note that $K$ depends on $m$ and $\lambda$.

\begin{Lem}
\label{LemK}
If $\lambda$ is a root of the characteristic polynomial~(\ref{polycar}) such that $\Re (\lambda)>-\frac 12$
and if $C$ is any complex number, then $K$ maps ${\cal M}_2(C) $ into itself.
\end{Lem}

\pff
Since $\Re (\lambda)>-1$, the random variable $e^{ -\lambda T}$ has an expectation.
Furthermore, by (\ref{beta}),  $m \g E e^{ -\lambda T} = 1$ as $\lambda$ is a root of (\ref{polycar}).
 This ensures the conservation of the expectation by $K$.
 Since $\Re (\lambda )>-\frac 12$, then $\g E |e^{ -\lambda T}|^2 < \infty$ and $K\mu$ admits a second absolute moment whenever
 $\mu$ does.    Therefore $K\mu \in   {\cal M}_2(C)$ whenever $\mu \in  {\cal M}_2(C). $
\QED

Notice that a solution of Eq.~(\ref{fixedpoint}) is a fixed point of $K$.
We shall use the Banach fixed point theorem for two different metrics on $\rond M_2(C)$ to study the existence
and uniqueness of solutions of Eq.~(\ref{fixedpoint}).

%%%%%%%%%%%%%%%%%%%%%%%%%%%%%
\subsection{Wassertein distance}
%%%%%%%%%%%%%%%%%%%%%%%%%%%%%%

Let $d_2$ be the Wasserstein distance on  ${\cal M}_2(C)$ (see for instance Dudley \cite{Dudley}):
for $  \mu, \nu \in {\cal M}_2( C)$,
\begin{equation*}
d_2(\mu, \nu) =\left( \min_{(X,Y)}\g E\left(\left| X-Y\right| ^2\right)\right)^{\frac 12},
\end{equation*}
where the minimum is taken over couples of random variables $(X,Y)$ having respective marginal distributions
$\mu$ and $\nu$;
the minimum is attained  by the Kantorovich-Rubinstein Theorem
-- see for instance Dudley~\cite{Dudley}, p.~421.
With this distance $d_2$, ${\cal M}_2(C)$ is a complete
metric space.

\begin{Th}
\label{W}
Let $\lambda\in\C$ be a root of the characteristic polynomial~(\ref{polycar}) such that $\Re(\lambda) >\frac 12$,
and let $C\in\C$.  Then $K$ is a  contraction on the complete metric space $({\cal M}_2(C), d_2)$, and the fixed
point equation (\ref{fixedpoint}) has a unique solution $Z$ in ${\cal M}_2( C)$.
\end{Th}

We now come back to the limit variable $W_1$ of $m$-ary search trees. Since  $\g E W_1= 1$ and $\g E|W_1|^2 <\infty$,
the following corollary is a direct consequence of Theorem~\ref{W}, applied for $\lambda =\lambda _2$.

\begin{Cor}
The distribution of the limit complex random variable $W_1$ is the unique solution in the space ${\cal M}_2(1)$ of
the fixed point equation (\ref{fixedpoint}).
\end{Cor}

\noindent
{\sc Proof of Theorem \ref{W}.}
We argue as in \cite{Gui90}, \cite{Ros92} and \cite{RR} where real random variables were considered.

By the Banach fixed point theorem, it suffices to show the contraction property.  Let $\mu, \nu \in {\cal M}_2( C)$.
Let $(X,Y)$ be a couple of complex-valued random variables such that $\rond L(X) = \mu$,
$\rond L(Y) = \nu$ and $d_2(\mu, \nu) = \sqrt{\g E |X-Y|^2}$.
Let  $(X_i,Y_i), i=1,\dots , m$ be $m$ independent copies of the $d_2$-optimal couple $(X,Y)$,
and $T$ be a real random variable with density $f_T$ defined by~(\ref{densitefT}), independent from all
$(X_i,Y_i)$.
Then,
$$
\rond L(e^{ -\lambda T } \sum_{i=1}^m X_i) = K\mu
{\rm ~~and~~~}
\rond L(e^{ -\lambda T } \sum_{i=1}^m Y_i ) =K\nu,
$$
so that
$$
\begin{array}{rcl}
d_2(K\mu, K\nu )^2
&\leq
&\displaystyle\g E\left|\left( e^{-\lambda T}\sum _{i=1}^mX_i\right)
-\left( e^{-\lambda T}\sum _{i=1}^mY_i\right)\right| ^2 \\ \\
&=&\g E\left| e^{ -\lambda T}\displaystyle\sum_{i=1}^m\left( X_i-Y_i\right)\right| ^2\\ \\
&=&\displaystyle\g E\left| e^{-\lambda T}\right| ^2 \g E\left|\sum_{i=1}^m\left( X_i-Y_i\right)\right| ^2\\ \\
&=&\displaystyle \g E\left| e^{ -\lambda T}\right| ^2\left(\sum _{i=1}^m \g E\left| X_i-Y_i\right| ^2
+\sum _{i\not= j}\g E\left( X_i-Y_i\right)\left( \overline{X_j}-\overline{Y_j}\right)\right)\\ \\
&=&m\g E\left| e^{ -2\lambda T}\right| \; d_2\left(\mu ,\nu\right) ^2.
\end{array}
$$
Since  $2\Re (\lambda )>1$, we have $m \g E\left| e^{ -2\lambda T}\right| <1$ (see (\ref{T})). Therefore $K$ is a contraction on $\rond M_2(C)$ and the proof is complete.

\QED

%%%%%%%%%%%%%%%%%%%%%%%%
\subsection{Distance defined with Fourier transforms}
%%%%%%%%%%%%%%%%%%%%%%%%%
We now give an alternative approach for the  characterization of the limit distribution via Fourier analysis.
We define another distance $d_2^*$ on $\rond M_2(C)$ as follows.
Take $\mu ,\nu\in\rond M_2(C)$ and denote respectively $\varphi$ and $\psi$ their characteristic functions.
By definition of $\rond M_2(C)$, both $\varphi$ and $\psi$ admit the expansion
$\varphi (t)=1+ i \langle t,C\rangle + O\left( |t|^2\right)$
when $t$ tends to $0$.
Therefore, one can define $d_2^*(\mu, \nu )$ by
$$
d_2^*(\mu ,\nu )= \sup_{t\in\g C\setminus\{ 0\}}  \frac{|\varphi(t) - \psi(t)| }{ |t|^2}.
$$
Clearly,  $d_2^* (\varphi, \psi) <\infty$, and $d_2^*$ is a distance on $\rond M_2 (C)$.
It can be easily checked that  $(\rond M_2 (C), d_2^*)$  is a complete metric space.

\vskip 5pt
The following result is a counterpart of  Theorem \ref{W}.
It gives an alternative proof for the existence and uniqueness of the solution of Eq. (\ref{fixedpoint}) in
the class of probability measures on $\g C$ with a given mean and finite second moments.

\begin{Th}
\label{contractiond2*}
Let $\lambda\in\C$ be a root of the characteristic polynomial~(\ref{polycar}) such that $\Re(\lambda) >\frac 12$,
and let $C\in\C$.  Then $K$ is a  contraction on the complete metric space $\left( \rond M_2(C), d_2^*\right)$,
and the fixed point equation (\ref{fixedpoint}) has a unique solution $Z$ in ${\cal M}_2( C)$.
\end{Th}

\pff
Thanks to Banach fixed point theorem, it suffices to prove that $K$ is a contraction on  ${\rond M_2 (C)}$ equipped with the metric $d_2^*$.
Let $\mu ,\nu\in\rond M_2(C)$ and let $\varphi$ and $\psi$ be their respective characteristic functions.
An elementary computation shows that the Fourier transform of $K\mu$ is
$t\mapsto \g E \varphi^m (e^{-\overline {\lambda} T} t)$ with a corresponding formula for $\nu$.
We have $|\varphi| \leq 1$,   $|\psi| \leq 1$, so that
\begin{equation*}
\g E |\varphi^m (e^{-\overline {\lambda} T} t) - \psi^m (e^{-\overline {\lambda} T} t)|
\leq
m \g E  |\varphi (e^{-\overline {\lambda} T} t) - \psi (e^{-\overline {\lambda} T} t)|.
\end{equation*}
Together with the inequality  $|\varphi (z) - \psi (z)| \leq  d_2^*  (\mu ,\nu ) |z|^2 $ applied to $z =  e^{-\overline {\lambda} T} t$,
this implies that
\begin{equation*}
d_2^* (K\mu , K\nu ) \; \leq \;  m \g E \left( e^{-2 \Re ({\lambda}) T}\right)  d_2^* (\mu , \nu ).
\end{equation*}
Since $2 \Re ({\lambda}) >1$, we have $m \g E\left( e^{-2 \Re ({\lambda}) T}\right) <1$ (see (\ref{T})). Therefore the above inequality shows that $K$ is a contraction on $\left( \rond M_2 (C),d_2^*\right)$.   \QED

\begin{Rem}
Denote $\rond F_2(C)$ the space of Fourier transforms of elements of $\rond M_2(C)$.
When $\lambda$ is a root of the characteristic polynomial~(\ref{polycar}) such that $\Re(\lambda) >\frac 12$,
the smoothing transformation $K$ can be identified as a map (also denoted by~$K$) on $\rond F_2(C)$
given by
\begin{equation}
\label{defKfunctions}
(K\varphi) (t) :=  \g E \varphi^m (e^{-\overline {\lambda} T} t),     \qquad t\in \C.
\end{equation}
The proof of Theorem~\ref{contractiond2*} also shows that $K$ is a contraction of $\rond F_2(C)$ for the
metric (also denoted by $d_2^*$) defined on $\rond F_2(C)$ by
\begin{equation}
\label{d2*Fonctions}
d_2^* (\varphi, \psi) := \sup_{t\neq 0}  \frac{|\varphi(t) - \psi(t)| }{ |t|^2}.
 \end{equation}

 Let $\rond D_2(C)$ be the space of all continuous functions $\varphi :\g C\to\g C$ that admit an expansion
 $\varphi (t)=1+ i \langle t,C\rangle + O\left( |t|^2\right)$ at $0$ and such that $||\varphi ||_\infty\leq 1$.
 Clearly, $\rond D_2(C)$ contains $\rond F_2(C)$.
 One can show that formula~(\ref{defKfunctions}) defines a mapping from $\rond D_2(C)$ into itself and
 that $K$ is a contraction for the metrics defined by~(\ref{d2*Fonctions}).
 This provides a proof of existence and unicity of solutions of~(\ref{fixedpointF}) on $\rond D_2(C)$.
\end{Rem}

\begin{Rem}
\label{remEqConvolution}
One can deal with the convolution equation~(\ref{psi}) by arguments of the same vein.
Similar computations show that this equation has a unique solution in the space $\rond E_2(\lambda ,C)$
of continuous functions $\Phi :\g C\to\g C$ that admit an expansion
 $\Phi (z)=1+ i\langle e^{{\lambda}z},C\rangle + O\left( e^{2{\lambda}z}\right)$ when $|z|$ tends to $+\infty$
 and such that $||\Phi ||_\infty\leq 1$.
This result appears one again as a consequence of Banach Theorem on $\rond E_2(\lambda ,C)$ for the metric
$$
d(\Phi, \Psi ) = \sup_{z\in\C}\left | \frac{\Phi (z) - \Psi (z)}{e^{2{\lambda}z}}\right |.
$$
As a consequence, this shows in particular that the Fourier (complex) transform $\varphi$  of $W_1$ satisfies:
for any $w\in\C^*$ and for any determination of the logarithm,
\begin{equation*}%\label{carW}
\varphi (w) =  \Phi \left(\frac {\log w}{\overline{\lambda_2}}   \right),
\end{equation*}
where $\Phi$ is the unique solution in $\rond E_2(\overline{\lambda _2},1)$ of Eq. (\ref{psi}).
This result is the reversed version of the change of variable $\Phi (z)=\varphi (e^{\overline{\lambda _2}z})$ that led from~(\ref{fixedpointmaryF}) to~(\ref{psi}).
\end{Rem}

%%%%%%%%%%%%%%%%%%%%%%%%%%%%%
\section{Density and support}
\label{secDensite}
%%%%%%%%%%%%%%%%%%%%%%%%%%%%%%

\renewcommand{\theequation}{\thesection.\arabic{equation}}
\setcounter{equation}{0}

In this section we prove results on the absolute continuity and on the support of solutions of the distributional equation
(\ref{fixedpoint}) via Fourier analysis. As applications, we show that the distribution of the  limit variable $W_1$ of the multitype branching process
$(X(t))$ has always a density and that its support is the whole complex plane.

\begin{Th}
\label{density} Let $\lambda$ be a complex number such that $\lambda\not= 1$  and  $\Re(\lambda) > 0$. Let
 $Z$ be a complex-valued random variable solution of the distributional equation (\ref{fixedpoint})
$$
Z \egalLoi e^{ -\lambda T } (Z^{(1)} + \cdots + Z^{(m)}),
$$
with  $\g E |Z|<\infty$ and $\g E Z\not= 0$.
Then the following assertions hold:
\begin{enumerate}
\item[(i)] the support of $Z$ is the whole complex plane $\C$;
\item[(ii)] as $|t| \rightarrow \infty$,  $\g Ee^{i\langle t,Z\rangle } = O(|t|^{-a})$, for any $a\in ]0, \displaystyle\frac1{\Re (\lambda)}[$;
\item[(iii)]
the distribution of $Z$ has a density with respect to the Lebesgue measure on~$\C$.
\end{enumerate}
\end{Th}
\begin{Rem}
\label{lambdaegale1}
When $\lambda = 1$, the distributional equation (\ref{fixedpoint}) becomes
\begin{equation}
\label{fixedpointlambdaegale1}
X \egalLoi e^{ -T } (X^{(1)} + \cdots + X^{(m)}).
\end{equation}
By Section \ref{smoothing}, it admits a unique solution in the space $\rond M_2(C)$ of probability measures on $\g C$, with a given mean $C$. Moreover, from the dislocation equations (\ref{dislocationfinitetime}), a similar argument shows that
$$
\xi:= \lim _{t\to +\infty}e^{-t}u_1\left( X(t)\right)
$$
is a solution of this equation. By Theorem \ref{thAsymptotiqueCT}, $\xi $ is Gamma-distributed. Therefore  the unique solution of (\ref{fixedpointlambdaegale1}) in $\rond M_2(C)$ is $C\gamma $ where $\gamma$ is Gamma$(1)$- distributed, and its support is the half line $C\g R_+$.
\end{Rem}

The following corollary gives the main result for the limit variable $W_1$ of the multitype branching process. It is a direct consequence of Theorem  \ref{density} since $\g EW_1 = 1$.

\begin{Cor}
\label{corDensite}
The distribution of $W_1$ admits a density with respect to the Lebesgue measure on  $\C$, and
its support  is the whole complex plane $\C$.
Moreover, as $|t| \rightarrow \infty$, $ \g Ee^{i\langle t,W_1\rangle } = O(|t|^{-a})$ for each $a\in ]0, \displaystyle\frac1{\Re (\lambda_2)}[$.
\end{Cor}

The proof of Theorem \ref{density} runs along the following lines. Let $\varphi$ be the Fourier transform of any solution $Z$ of (\ref{fixedpoint}). We prove that $\varphi$ is in $L^2(\g C)$  because it is dominated by $|t|^{-\delta}$ for some $\delta >1$ so that the inverse Fourier-Plancherel transform provides a square integrable density for $Z$. The guiding idea consists in adapting usual methods (developed in \cite{Liu99} and \cite{Liu01}) used for positive real-valued random variables to our complex-valued case thanks to the function defined for $r\geq 0$ by
\begin{equation*}
\psi(r) = \max_{|t| = r}|\varphi(t)|.
\end{equation*}
From now on,
$$A=e^{-\lambda T}.$$

We proceed by a series of lemmas.

 The first lemma concerns a property of the support of $Z$. For a complex-valued random variable $Z$ and for a complex number $z$, by definition,
$$
z\in \Supp (Z) \Longleftrightarrow \forall\varepsilon >0, \g P (|Z-z|\leq \varepsilon ) >0.
$$
\begin{Lem}
\label{lemma0}
Let $z\in\C$. Then
$$z\in \Supp (Z) \Longrightarrow D(0, |z|)\subseteq \Supp (Z),$$
where $D(0, |z|)$ denotes the open disc with center $0$ and radius $|z|$.
\end{Lem}

\noindent
{\sc Proof.}
We first prove the following implication:
%$\forall a>0$, $\forall w\in\g C$,
\begin{equation*}
\big[ a\in\Supp (A)~{\rm and}~z\in\Supp (Z)\big] \Longrightarrow maz\in\Supp (Z).
\end{equation*}
Indeed, let $\varepsilon >0$, $a\in\Supp (A)$ and $z\in\Supp (Z)$.
Let also $Z^{(1)},\dots ,Z^{(m)}$ be independent copies of $Z$.
Then, with positive probability, $|A-a|\leq\varepsilon$ and $|Z^{(k)}-z|\leq\varepsilon$ for any $k$.
Therefore, with positive probability,
$$
\begin{array}{rl}
\left| A(Z^{(1)}+\dots +Z^{(m)})-maz\right|
&\displaystyle=\left| mz(A-a)+A\sum _{k=1}^m(Z^{(k)}-z)\right|\\ \\
&\displaystyle\leq m\varepsilon |z|+\left( |a|+\varepsilon\right)m\varepsilon .
\end{array}
$$
The positive real $\varepsilon$ being arbitrary, this shows that $ma z\in\Supp A(Z^{(1)}+\dots +Z^{(m)})$ which
implies that $maz\in\Supp (Z)$ by Eq.~(\ref{fixedpoint}), proving the claim.

Let $z\in\Supp (Z)$.
%$z\neq 0$.
Since $\Supp (T)=\g R_+$ (see~(\ref{densitefT})), the claim implies that for any $t\geq 0$,
$me^{-\lambda t}z\in\Supp (Z)$.
Iterating this property of $\Supp (Z)$ shows that
\begin{equation}
\label{grosSupport}
\left\{ m^ne^{-\lambda t}z,~n\in\g N,~t\in\g R_+\right\}\subseteq\Supp (Z).
\end{equation}
Since the support of a probability measure is a closed set, to show that
$D(0,|z|)\subseteq\Supp (Z)$ it suffices to prove that
$\left\{ m^ne^{-\lambda t},~n\in\g N,~t\in\g R_+\right\}$ is everywhere dense in the unit disc.
Taking  logarithm,
we show hereunder that
$$
\rond G:=\g N\log m+2i\pi\g N-\lambda\g R_+
$$
is everywhere dense in the half-strip
$$
\rond B:=\left\{ x+iy,~x<0,~-2\pi < y\leq 0\right\}
$$
which implies the desired  result.

Let $\sigma$ and $\tau$ denote respectively the real and imaginary parts of $\lambda$.
Remember that, as soon as $Z$ is a solution of Eq.(\ref{fixedpoint}) such that $\g E |Z|<\infty$ and $\g E Z\not= 0$, then $\lambda$ is a root of $\chi _{R_G}$ defined by (\ref{polycar}); this implies that $\lambda$ is an algebraic number, and so is $\sigma /\tau$.

Let us prove that  $\rho:=i\pi/\log m$ is a transcendental number. In fact, if $\rho$ were algebraic, then  by the Gelfond-Schneider theorem\footnote{The Gelfond-Schneider theorem states that if $a$ and $b$ are algebraic numbers with $a\neq 0,1$ and if $b$ is not a rational number, then any value of $a^b= \exp (b \log a) $ is a transcendental number.},
  $m^\rho$ would be a transcendental number; but this is impossible because  $m^\rho= \exp (i\pi) =  -1$.
Therefore $\rho$ is a transcendental number.

It follows that  $2\pi\sigma /(\tau\log m) $  is not an algebraic number, hence not a rational number.
So  by a classical result,
$$
\g N\log m-\g N\frac{2\pi\sigma}{\tau}  \;\;  \mbox{ is a dense subset of }\;  \g R.
$$
Let $b=x+iy\in\rond B$ with $x<0$ and $-2\pi<y\leq 0$,  and let $\varepsilon >0$.
Approximating the real number $x-\frac {y\sigma}{\tau}$ by an element of
$\g N\log m-\g N\frac{2\pi\sigma}{\tau}$, we take $n,k\in\g N$ such that
$$
\left| \big( n\log m-k\frac{2\pi\sigma}{\tau}\big) -\big( x-\frac {y\sigma}{\tau}\big)\right|\leq\varepsilon.
$$
Therefore
$$ |b-g|\leq\varepsilon, \mbox{ where }
g=n\log m+2ik\pi-t\lambda\in\rond G \mbox{ with } t=\frac{2k\pi -y}{\tau}\geq 0.
$$
This completes the proof of Lemma~\ref{lemma0}.
\QED
Lemma \ref{lemma0} leads to the following key property of $\psi$,
which will imply that  the characteristic function $\varphi$  of $Z$ satisfies the Cram\'er condition and then the Riemann-Lebesgue condition.
\begin{Lem}
\label{lemma1}
$\forall r>0, \psi(r) <1 $ .
\end{Lem}
{\sc Proof.}
Obviously, $\psi (0)=1$ and $\psi (r)\leq 1$ for any $r\geq 0$.
Suppose that $r_0>0$ is such that $\psi (r_0)=1$.
Take thus $z_0\in\g C$ and $\theta _0\in\g R$ such that
$$
|z_0|=r_0~~{\rm and}~~\g Ee^{i\langle z_0,Z\rangle}=e^{i\theta _0}.
$$
The complex random variable $e^{i(\langle z_0,Z\rangle -\theta _0)}$ is of mean $1$ and takes its values on the unit disc, so
that it is almost surely equal to $1$.
This implies that $\Supp (Z)$ is contained in a set of countably many  parallel lines of the complex
plane.
This contradicts Lemma~\ref{lemma0} since such a set of lines is negligible
with respect to the Lebesgue measure on $\g C$.
\QED
\begin{Rem}
The preceding arguments show the following assertion: for any complex-valued random variable  $Z$, if   $| \g E(e^{i\langle z_0,Z\rangle}) |=1$ for some  $z_0\in\g C\setminus\{ 0\}$,  then $\Supp (Z)\subseteq a+b\g Z+c\g R$ for some  $a,b,c\in\g C$ (a set of countably many  parallel lines).
The algebraicity of $\lambda$ that leads to the proof of Lemma~\ref{lemma0} can thus be seen as a nonlattice assumption on the fixed point equation~(\ref{fixedpoint}).
\end{Rem}
We now prove that the characteristic function $\varphi$ satisfies the Riemann-Lesbesgue condition.
\begin{Lem}
\label{lemma2}
$\displaystyle\lim_{r\rightarrow +\infty} \psi(r) = 0$.
\end{Lem}
\noindent {\sc Proof.} We argue as in the proof of Theorem 3.1 of \cite{Liu99} or
Lemma 3.1 of \cite{Liu01}.  Notice that from the distributional equation (\ref{fixedpointF}) we have
\begin{equation}
\label{fixedpointpF}
\psi(r) \leq \g E(\psi^m(r|A|)).
\end{equation}
$\bullet$ We first prove that $\displaystyle\limsup_{r\to +\infty} \psi(r) = 0$ or $1$. By Fatou's lemma,
$$
\limsup_{r\to +\infty} \psi(r) \leq \g E \limsup_{r\to +\infty} \psi^m(r|A|) = \left( \limsup_{r\to +\infty}\psi(r) \right) ^m,
$$
the last equality coming from $\g P(|A| >0) = 1$.
So the real number $l:=\limsup_r \psi(r)$ satisfies both $l\leq 1$ and $l\leq l^m$; this implies that $l=0$ or $1$.

\medskip
\noindent $\bullet$
Suppose that $\limsup_r \psi(r) =1$. By Lemma \ref{lemma1}, $\psi(1)<1$. For any $\varepsilon \in ] 0, 1-\psi(1)[$, define
$$
\left\{
\begin{array}{l}
r_1(\varepsilon) = \max \{ r\in]0,1[, \psi(r) = 1-\varepsilon\},  \\ \\
r_2(\varepsilon) = \min \{ r>1, \psi(r) = 1-\varepsilon\} .
\end{array}
\right.
$$
These quantities are well defined because $\psi(0) = 1$ and $\psi$ is continuous.
Then $\psi(r_1(\varepsilon)) = \psi(r_2(\varepsilon)) = 1-\varepsilon  $ and
for any $ r\in[r_1(\varepsilon),r_2(\varepsilon)], \psi(r) \leq 1-\varepsilon.$

Let us prove that $r_1(\varepsilon)$ goes to $0$ when $\varepsilon$ tends to $0$.
Take any limit point $\rho$ of $r_1(\varepsilon)$. Since $\psi$ is continuous,
$\psi(\rho) = 1-\varepsilon$ which implies by Lemma \ref{lemma1} that $\rho = 0$:
the only possible limit point is $0$.

By (\ref{fixedpointpF}), we have
$$
\psi(r) \leq \g E \psi(r|A|).
$$
Iterating this inequality we see   that for all $n\geq 1$,
$$
\psi(r) \leq \g E \psi(r|A_1|\dots |A_n|),
$$
where $(|A_i|)_{i\geq 1}$ are independent copies of $|A|$.
With the notation
$$
\lambda_n(r,\varepsilon):= \g P\Big( r_1(\varepsilon) < r |A_1|\dots |A_n| \leq r_2(\varepsilon)\Big),
$$
we have for any $r>0$
$$
\psi(r) \leq (1-\varepsilon) \lambda_n(r,\varepsilon) + 1- \lambda_n(r,\varepsilon)
= 1 - \varepsilon \lambda_n(r,\varepsilon).
$$
Again by  (\ref{fixedpointF}),
$$
1-\varepsilon =  \psi(r_2(\varepsilon)) \leq \g E \psi^m\Big( r_2(\varepsilon) |A|\Big) \leq \g E \left( 1 - \varepsilon \lambda_n\Big(r_2(\varepsilon)|A|,\varepsilon\Big)\right)^m.
$$
In other words
\begin{equation}
\label{star17juillet}
\frac{\g E \left( 1- \left( 1 - \varepsilon \lambda_n\Big(r_2(\varepsilon)|A|,\varepsilon\Big) \right)^m\right)}{\varepsilon}
\leq 1.
\end{equation}
We are going to pass to the limit in the above ratio when $\varepsilon$ tends to $0$. Rewrite
$$
\lambda_n\Big(r_2(\varepsilon)|A|,\varepsilon\Big)  = \g P \left( \frac{r_1(\varepsilon)}{r_2(\varepsilon)} < |A||A_1| \dots |A_n| \leq 1\right),
$$
and remember that $r_1(\varepsilon)\leq 1\leq r_2(\varepsilon)$ and that $r_1(\varepsilon)$ goes to $0$ when $\varepsilon$ tends to $0$, so that $\frac{r_1(\varepsilon)}{r_2(\varepsilon)} \leq \frac{r_1(\varepsilon)}1$ goes to $0$ when $\varepsilon$ tends to $0$. Consequently,
$$
\lambda_n\Big(r_2(\varepsilon)|A|,\varepsilon\Big) \  \convepsilon \ \g P\Big( 0\leq  |A||A_1| \dots |A_n| \leq 1\Big) = \mu_n(|A|) \hskip 5mm a.s.,
$$
where, for any $x>0$,
$$
\mu_n(x) := \g P\Big( x|A_1| \dots |A_n| \leq 1\Big) .
$$
Therefore
$$
\frac{1- \left( 1 - \varepsilon \lambda_n\Big(r_2(\varepsilon)|A|,\varepsilon\Big) \right)^m}{\varepsilon} \  \convepsilon \ m\mu_n(|A|) \hskip 5mm a.s..
$$
The above ratio is a function of $\varepsilon$, uniformly bounded on the compact set $[0, 1-\psi(1)]$, so that by dominated convergence
and (\ref{star17juillet}),
\begin{equation}
\label{majorationm}
\frac{\g E \left( 1- \left( 1 - \varepsilon \lambda_n\Big(r_2(\varepsilon)|A|,\varepsilon\Big) \right)^m\right)}{\varepsilon}  \  \convepsilon \ m\g E\mu_n(|A|)\leq 1.
\end{equation}
Besides, by Markov inequality
$$
1-\mu_n(x)\leq x\g E(|A_1| \dots |A_n|) = x\left( \g E |A| \right)^n.
$$
Since $\Re(\lambda) >0$,  $\g E |A| = \g E | e^{-\lambda T}| <1 $ (see (\ref{beta})), which implies that $\lim_{n\rightarrow \infty} \mu_n(x) = 1$,
so that
$$\lim_{n\rightarrow \infty} \g E\mu_n( |A| ) = 1
$$
 by dominated convergence.
  This contradicts (\ref{majorationm}) because $m\geq 2$.
\QED

\bigskip
We shall need an information  about the decay rate of $\varphi (t)$, of the form  $ \varphi (t) = O (|t|^{-\delta})$ for some $\delta > 0$ when $|t| \rightarrow \infty$. To this end, we shall use the following Gronwall-type technical Lemma of \cite{Liu99} (see also  Lemma 3.2 in \cite{Liu01}).
\begin{Lem}
\label{lemma3} \cite{Liu99}
Let $\psi: \R_+ \rightarrow \R_+$ be a  bounded function
and let $B$ be a positive random variable such that for some constants
$p \in ]0,1[, a >0, C \geq 0$ and for all $r>0$,
$$ \psi (r) \leq p\g E\psi (Br) + Cr^{-a}.  $$
 If $p\g E(B^{-a})<1$, then $\psi (r) = O(r^{-a})$ as $r\rightarrow \infty$.
\end{Lem}

\noindent
%{\sc Proof of Lemma \ref{lemma3}.}
This is  Lemma 4.1 of   \cite{Liu99}. It can be proved as follows.
%%%% pour Q : on a simplifiŽ
%In fact,  we can assume $r_0=0$ by taking $C$ large enough.
Let $\{B_i\} $ be independent copies of $B$. Then
by induction, for all $n\geq 1$ and all $r>0$,
$$ \psi (r) \leq p^n \g E\psi (B_1...B_n r) + Cr^{-a}
                   [1+ p\g E(B^{-a}) +... + ( p\g E(B^{-a}))^{n-1}].  $$
                    Letting $n\rightarrow \infty$
we see that  for all $r>0$,
$$ \psi (r) \leq Cr^{-a} / [1-p\g E(B^{-a})]. $$

\begin{Lem}
\label{lemma4}
For  all $a\in]0,\displaystyle\frac 1{\Re(\lambda)}[$,  as $r \rightarrow \infty$,
\begin{equation*}
%\label{eqlemma4}
\psi(r)  =  O (r^{-a}) .
\end{equation*}
\end{Lem}

\pff
We have already seen from the distributional equation (\ref{fixedpointF}) that
\begin{equation*}
%\label{fixedpointpF}
\psi(r) \leq \g E(\psi^m(r|A|)),
\end{equation*}
where $A=e^{ -\lambda T }$.
By Lemma \ref{lemma2}, for any $\varepsilon >0$, there is some $r_{\varepsilon} >0$
such that $\forall r \geq r_{\varepsilon}$, $\psi(r) \leq \varepsilon$. So
$$
\psi(r)\leq \g  \varepsilon^{m-1} \g E \psi(r|A|) + \g P(r|A| \leq r_{\epsilon}).
$$
Therefore by Markov inequality, for  $a\in]0,\displaystyle\frac 1{\Re(\lambda)}[$,
$$
\psi(r)\leq   \varepsilon^{m-1} \g E \psi(r|A|) +  r^{-a} (r_{\varepsilon})^a \g E(|A|^{-a}).
$$
By (\ref{beta}), $\g E(|A|^{-a}) = (m-1)B(1- a\Re(\lambda), m-1) < \infty$.
Taking $\varepsilon >0$ small enough such that $ \varepsilon^{m-1} \g E(|A|^{-a})<1$, we see that the desired result follows from Lemma \ref{lemma3} and the
preceding inequality on $\psi(r)$.~\QED

\bigskip
We can now finish the proof of Theorem \ref{density}.

\medskip
\noindent {\sc Proof of Theorem \ref{density}}.
Part (i) of the theorem comes from two facts as shown in the following.

On the one hand, by Lemma \ref{lemma0}, as soon as $z\in\C$ is a point in the support of $Z$, we have
$D(0, |z|)\subseteq \Supp (Z)$, where $D(0, |z|)$ denotes the open disc with center $0$ and radius $|z|$.

 On the other hand, the support of $Z$ is unbounded. Indeed, as in (\ref{grosSupport}), at the beginning of the proof of Lemma \ref{lemma0}, as soon as $z\in\C\setminus \{0\}$ is a point in the support of $Z$, for any $t>0$ and for any $n\in\N$, $m^n e^{-\lambda t}z$ is in the support of $Z$.

For Part (ii), notice that by Lemma \ref{lemma4},
for  all $a\in]0,\displaystyle\frac 1{\Re(\lambda)}[$,
\begin{equation}
\label{eqlemma4bis}
\varphi(t)  =  O (t^{-a})    \qquad \mbox{ as } |t| \rightarrow \infty.
\end{equation}
Since $\g E Z\not= 0$, by Eq. (\ref{fixedpointF}) we obtain  $m \g Ee^{-\lambda T}=1$, hence $m\g E  | e^{-\lambda T} | = m\g Ee^{-\Re (\lambda) T}  > 1 $ as soon as $\Im(\lambda)\not= 0$. Notice that if $\Im(\lambda) = 0$, then $\lambda = 1$ by the equation $m \g Ee^{-\lambda T}=1$. So the hypotheses $\lambda\not= 1$
and $\g E Z\not= 0$ imply that $\Im(\lambda)\not= 0$ and $\Re( \lambda) <1$ (cf. (\ref{T})).  It follows that
   (\ref{eqlemma4bis}) holds for some   $a>1$, so that the Fourier transform $\varphi$ of $Z$ is in $L^2$.
 Therefore by the inversion formula of Fourier-Plancherel transform, the distribution of $Z$  has a density in  $L^2$  with respect to the Lebesgue measure on $\C$.
  This ends  the proof of  Theorem \ref{density}.
\QED

\begin{Rem}
\label{RemDensity}
In fact we have the following more general result.   Let $\lambda$ be an complex number  with  $\sigma:=\Re(\lambda) > 0$, $\tau:=\Im (\lambda) \not= 0$ and satisfying the arithmetical condition:
$$
\frac{\pi\sigma}{\tau\log m} \notin \g Q
$$
 and let $Z$ be a nontrivial solution of Eq. (\ref{fixedpoint}) (with or without first moment).
Then the distribution of $Z$ is absolutely continuous with respect to the Lebesgue measure on $\C$, and its support is the whole complex plane $\C$.
\end{Rem}
\noindent
%{\sc Proof of Remark \ref{RemDensity}.}
 To see the conclusions of Remark \ref{RemDensity}, we can argue as follows.
 In the general case where the expectation of $ Z$ may not exist, Lemma \ref{lemma0} still holds thanks to the arithmetical condition. The remaining of the proof is the same, except at the end, where $\Re(\lambda)>1$ is no more ensured. Nevertheless we have an additional argument by iteration.
 Iterating the distributional equation (\ref{fixedpoint}), we obtain for $n\geq 1$,
\begin{equation*}
%\label{fixedpointIter1}
Z   \egalLoi \sum_{u_1\dots u_n \in \{1,\dots , m\}^n} A A_{u_1} \dots A_{u_1 \dots u_{n-1}} Z^{(u_1\dots u_n)},
\end{equation*}
where $A= e^{-\lambda T}$, $A_u$ are independent copies of $A$ (indexed by finite sequences of integers $u$),
$Z^{(u)} $ are independent copies of $Z$, the two families $\{A_u\}$ and $\{Z^{(u)}\}$ are also independent of each other;
by convention,
$A_{u_1} \dots A_{u_1 \dots u_{n-1}}$ is taken to be $1$ when $n=1$.  It is convenient to rewrite this equation in the form
\begin{equation}
\label{fixedpointIter2}
      Z   \egalLoi \sum_{j=1}^{m^n} Y_j Z^{(j)} ,
\end{equation}
where $Z^{(j)}$ are independent copies of $Z$ which are also independent of $\{Y_j\}$.
For fixed $y= (y_j: 1\leq j \leq m^n)$ with $\prod_{j=1}^{m^n}y_j \not= 0$, by Lemma \ref{lemma4}, for $a \in ]0, 1/\Re (\lambda)[$ and some constant $c>0$,
\begin{eqnarray*}
  \left| \g E \exp { (i \langle t,  \sum_{j=1}^{m^n} y_j Z^{(j)} \rangle) }\right|
     &\leq &   \prod_{j=1}^{m^n} c |t y_j|^{-a}  \\
      &=& C(y) |t|^{-m^n a},
\end{eqnarray*}
where $C(y) =\prod_{j=1}^{m^n}c |y_j|^{-a} >0$ does not depend on $t$.
Let $n\geq 1$  be large enough such that $m^n a >1$. Then  $\sum_{j=1}^{m^n} y_j Z^{(j)}$ is absolute continuous (with respect to the Lebesgue measure on $\C$) as its Fourier transform is square integrable on $\C$. This implies that for each Borel set $B$ of $\C$ with Lebesgue measure $0$, we have
\begin{equation*}
      \g P ( \sum_{j=1}^{m^n} y_j Z^{(j)} \in B) =0.
\end{equation*}
It follows from  Eq. \ref{fixedpointIter2} (by conditioning on $(Y_j)$) that $\g P(Z \in B) = 0$.

%%%%%%%%%%%%%%%%%%%%%%%%%%%%%%%%
\section{Exponential moments and Laplace series}
\label{laplace}
%%%%%%%%%%%%%%%%%%%%%%%%%%%%%%%%

\renewcommand{\theequation}{\thesection.\arabic{equation}}
\setcounter{equation}{0}

In this section we consider a solution $Z$ of Eq. (\ref{fixedpoint}) and we show that its exponential moments exist in a neighborhood of $0$, so that the moment exponential generating series of $Z$ defines an analytic function in a neighbourhood of the origin.
We show that this function satisfies a very simple differential equation.

\begin{Th}
\label{Expmoments}
Let $\lambda \in \g C$ be a root of the characteristic polynomial~(\ref{polycar}) with $\Re (\lambda) >1/2$ and let $Z$ be a solution of Eq. (\ref{fixedpoint}).
There exist some constants $C>0$ and $\varepsilon >0$
such that for all $t\in \C$ with $ |t| \leq \varepsilon$,
 \begin{equation}
 \label{bdExpmoments}
  \g Ee^{ \langle t,Z\rangle   } \leq      e^{ \Re (t)  + C |t|^2 }  \; \mbox{ and } \;
  \g Ee^{ |tZ | } \leq   4   e^{ |t|  +  2C |t|^2 } .
 \end{equation}
 \end{Th}

To prove this theorem, we use Mandelbrot's cascades in the complex setting (see Barral et al. \cite{Ba10} for independent interest about complex Mandelbrot's cascades).
We still denote  $A= e^{-\lambda T}$. Then $m\g EA=1$ because $\lambda $ is a root of the characteristic polynomial~(\ref{polycar}) and $m\g E|A|^2 < 1$ because $\Re (\lambda) >1/2$ (see (\ref{T})).
Let $A_u, u\in U$ be independent copies of $A$, indexed by all finite sequences of integers
$$
u = u_1...u_n \in U:= \bigcup_{n\geq 1} \{ 1, 2, \dots , m\}^n
$$
and set $Y_0=1$, $Y_1 = mA$ and for $n\geq 2$,
 \begin{equation*}
 Y_n =  \sum_{u_1...u_{n-1} \in \{1,..., m\}^{n-1}}  mAA_{u_1} A_{u_1u_2} \dots A_{u_1 ... u_{n-1}}.
 \end{equation*}
As $m\g EA=1$, $(Y_n)_n$ is a martingale.
This martingale has been studied by many authors in the real random variable case, especially in the context of Mandelbrot's cascades,
see for example $\cite{Liu01}$  and the references therein.
It can be easily seen that
\begin{equation}
\label{EqYn}
    Y_{n+1} = A \sum_{i=1}^m Y_{n,i}
\end{equation}
where $ Y_{n,i} $ for $1\leq i \leq m $ are independent of each other and independent of $A$ and  each has the same distribution as $Y_n$.
Therefore for $n\geq 1$, $Y_n$ is square-integrable and
$$ \Var Y_{n+1} = (\g E |A|^2 m^2 - 1) + m\g E|A|^2 \Var Y_n,  $$
where $\Var X = \g E\left( |X-\g EX|^2\right)$  denotes the variance of $X$.
Since $m\g E|A|^2 < 1$, the martingale  $(Y_n)_n$ is bounded  in $L^2$, so that the following result holds.

\begin{Lem}
\label{YnL2}
Let $\lambda $ be a root of the characteristic polynomial~(\ref{polycar})  with $\Re (\lambda) >1/2$.  Then, when $n\rightarrow + \infty$,
 \begin{equation*}
    Y_n \rightarrow  Y_\infty \mbox{ a.s. and in } L^2,
  \end{equation*}
where $Y_\infty$ is a (complex-valued) random variable with variance
 \begin{equation*}
  \Var (Y_\infty) = \frac {\g E |A|^2 m^2 - 1}{1  - m\g E|A|^2 }.
 \end{equation*}
\end{Lem}

Notice that, passing to the limit in~(\ref{EqYn}) gives a new proof of the existence of a solution $Z$ of  Eq. (\ref{fixedpoint}) with $\g EZ=1$ and finite second moment whenever $\Re (\lambda) >1/2$. From Section \ref{smoothing}, we have the unicity of solution of this equation so that Theorem~\ref{Expmoments} is proved as soon as it holds for $Y_\infty$.

\begin{Lem}
\label{expmomentsY}
Under the condition of Lemma  \ref{YnL2},  there exist some constants $C>0$ and $\varepsilon >0$ such that for
all $t \in \C$ with $|t| \leq \varepsilon $, we have
\begin{equation}
 \label{bound-phi-inf}
  \g E  e^{\langle t, Y_\infty \rangle  } \leq e^{\Re (t)  + C |t|^2 }.
 \end{equation}
\end{Lem}

\pff
As in  \cite{Ros92} and \cite{LR00} (where a similar problem for real random variables  was considered), we use an induction argument. Notice that by Eq. (\ref{EqYn}),
writing
$$ \varphi_n (t) := \g E e^{ \langle t, Y_{n} \rangle },  \qquad t \in \C, n \geq 0,  $$
we have
\begin{equation}
\label{EqPhi-n}
 \varphi_{n+1} (t) =  \g E \varphi_n^m (\overline A t), \qquad t \in \C.
\end{equation}
We shall prove that there exist some constants $C>0$ and $\varepsilon >0$ such that for all $n \geq 0$ and
all $t \in \C$ with $|t| \leq \varepsilon $, we have
\begin{equation}
 \label{bound-phi-n}
 \varphi_n (t) \leq e^{\Re (t) + C |t|^2 } .
 \end{equation}
Let us prove (\ref{bound-phi-n})  by induction. The inequality holds clearly for $n=0$ since $\varphi_0 (t) = e^{\Re (t)}$.
Assume that it holds for some $n\geq 0$ and all $t \in \C$ with $|t|\leq \varepsilon$.
Then  writing $A= A_1 + iA_2$ $(A_i \in \R) $, using $|A | \leq 1 $ and Eq. (\ref{EqPhi-n}),
we have for  $t = t_1 + it_2 $ $(t_i\in \R)$ with
$|t| \leq \varepsilon $,
  \begin{eqnarray}
\label{recuBoundPhi-n}
 \varphi_{n+1} (t) & \leq &    \g E \exp \{ m (A_1 t_1 + A_2t_2 + C |A|^2 (t_1^2 + t_2^2)) \} \nonumber \\
                & =&  e^{\Re( t) + C |t|^2 } g(t_1,t_2),
\end{eqnarray}
where $ g(t_1,t_2) = \g Ee^{h(t_1,t_2)}$  with
$$  h(t_1,t_2) =  (m A_1-1) t_1 + mA_2t_2 + C(m|A|^2 -1) (t_1^2 + t_2^2) . $$
Notice that $g(0,0) = 1$. It remains to prove that (0,0) is a local maximum of $g$.
Clearly,
\begin{eqnarray*}
\frac{\partial g}{\partial t_i} &=& \g E e^h \left[   \frac{\partial h}{\partial t_i}\right],  \qquad i=1,2 \\
 \frac{\partial^2 g}{\partial t_i} &=& \g E e^h \left[   \big( \frac{\partial h}{\partial t_i}\big)^2 + \frac{\partial^2 h}{\partial t_i}\right] ,  \qquad i=1,2 \\
\frac{\partial^2 g}{\partial t_1 \partial t_2} &=& \g E e^h \left[  \frac{\partial h}{\partial t_1} \frac{\partial h}{\partial t_1} + \frac{\partial^2 h}{\partial t_1 \partial t_2}\right].
 \end{eqnarray*}
Notice that, a.s.
$$ \frac{\partial h}{\partial t_1} (0,0) =  (m A_1-1), \hskip 9mm  \frac{\partial h}{\partial t_2} (0,0) = mA_2, $$
$$\frac{\partial^2 h}{\partial t_1 \partial t_2} (0,0) =0,
\hskip 9mm \frac{\partial^2 h}{\partial t_i} (0,0) = 2C (m|A|^2 -1),  \qquad i=1,2.
$$
Recall that $m\g EA=1$, so that $m\g EA_1 = 1$ and $m\g EA_2 =0$; hence
$$ \frac{\partial g}{\partial t_1}(0,0) = \g E  (m A_1-1)=0, \hskip 9mm \frac{\partial g}{\partial t_2}(0,0) = \g E  (m A_2) =0, $$
so that $(0,0)$ is a critical point of $g$.
Moreover,
\begin{eqnarray*}
 \frac{\partial^2 g}{\partial t_1} (0,0) &=& \g E \left[ (m A_1-1)^2 + 2C (m|A|^2 -1)\right],  \\
 \frac{\partial^2 g}{\partial t_2} (0,0)&= &\g E \left[ (m A_2)^2 + 2C (m|A|^2 -1)\right] ,\\
 \frac{\partial^2 g}{\partial t_1 \partial t_2} (0,0) &=& \g E (m A_1-1) (mA_2).
\end{eqnarray*}
As $\g E(m|A|^2 -1) <0$ (recall that $\Re (\lambda) >1/2$),
it follows that  the Hessian matrix  at $(0,0)$ is definite negative for $C >0$ large enough which implies that $g(0,0)$ is a local maximum of $g$. So for $\varepsilon >0$ small enough, $g(t_1,t_2) \leq g(0,0) = 1$
for all $t = t_1 +i t_2 $ with $|t| \leq \varepsilon$. Hence by (\ref{recuBoundPhi-n}), for such $\varepsilon$ and $C$ which do not depend on $n$, (\ref{bound-phi-n}) holds  for $n+1$. Therefore, by induction,  it holds for all $n\geq 0$.

Letting $n\rightarrow \infty$ in   $ (\ref{bound-phi-n})$, we see that inequality (\ref{bound-phi-inf}) holds by Fatou's lemma. \QED

\medskip
\noindent{\sc Proof of Theorem \ref{Expmoments}.}  By the unicity of solution of Eq. (\ref{fixedpoint}),  $\rond L(Z)= \rond L(Y_\infty )$.
So by Lemma \ref{expmomentsY}, there are some constants
 $C>0$ and $\varepsilon >0$ such that the first inequality of (\ref{bdExpmoments}) holds. To show the second one,
 notice that   $|t| | \Re (Z) |  + |t| |\Im (Z)|$ takes one of the four values $\pm |t| \Re (Z)   \pm |t|  \Im (Z)$
 (according to the signs of
$\Re (Z) $ and $ \Im (Z)$), so that a.s.
 \begin{eqnarray*}
e^ { | t  Z | }   &\leq &  e ^ {|t| | \Re (Z) |  + |t| |\Im (Z)|} \\
                                  & \leq &  e ^ {| t|  \Re (Z)   + |t| \Im (Z) }
                                   +  e ^ {|t|  \Re (Z)   -|t|  \Im (Z) }
                                   + e ^ { -|t|  \Re (Z)   +|t|  \Im (Z) }
                                    + e ^ { - |t|  \Re (Z)   - |t| \Im (Z) }.
\end{eqnarray*}
Taking expectation in both sides, and noticing that $\scriptstyle\pm |t| \Re (Z)   \pm |t|  \Im (Z) = \langle (\pm 1\pm i)|t|, Z\rangle $, we see that the second inequality in (\ref{bdExpmoments}) follows from the first one. \QED
\vskip 10pt
Suppose that $Z$ is any solution of Eq. (\ref{fixedpoint}) under the assumptions of
Theorem~\ref{Expmoments}.
The second inequality~(\ref{bdExpmoments}) shows that the exponential generating series of absolute moments
of $Z$ has a positive radius of convergence so that the formal Laplace series
$$
L(z):=\sum _{p\geq 0}\frac{\g EZ^p}{p!}z^p
$$
defines an analytic function in a neighbourhood of the origin.
One can also write $L(z)=\g Ee^{zZ}$ when $|z|$ is sufficiently small.

Let's come back to the dislocation equations~(\ref{system}) satisfied by the limit variables $W_1,\dots ,W_m$.
These variables admit finite (absolute) moments at any order.
For any $k\in\{ 1,\dots ,m\}$, let $L_k$ be the formal Laplace series defined by
$$
L_k(z):=\sum _{p\geq 0}\frac{\g E\left( W_k^p\right)}{p!}z^p.
$$
The dislocation equations~(\ref{system}) imply recursive relations on $W_k$'s moments.
Developing these relations with the multinomial formula implies that  $L_k$ satisfy the formal differential
system
\begin{equation}
\label{systemeFormel}
\left\{
\begin{array}{l}
\displaystyle
\forall k\in\{1,\dots ,m-2\},~
L_k(z)+\frac{\lambda _2}kzL_k'(z)=L_{k+1}(z),
\\
\displaystyle
L_{m-1}(z)+\frac{\lambda _2}{m-1}zL_{m-1}'(z)=\left( L_{1}(z)\right) ^m ,
\end{array}
\right.
\end{equation}
with boundary conditions
\begin{equation*}
%\label{bordFormel}
\left\{
\begin{array}{l}
\displaystyle
L_k(0)=1,~1\leq k\leq m-1 , \\
\displaystyle
L_k'(0)=\g E(W_k)=u_2\left( X_k(0)\right)=\binom{\lambda _2+k-1}{k-1}.
\end{array}
\right.
\end{equation*}
Since $W_1$ satisfies the assumptions of Theorem~\ref{Expmoments}, the series  $L_1$ has a positive radius
of convergence as shown above.
Therefore, the same holds for all $L_k$ because of the system~(\ref{systemeFormel}) so that the $L_k$ define,
near the origin, analytic functions related by (\ref{systemeFormel}).

Let $\rho$ be any complex $(m-1)$-th root of $(-1)^m(m-1)!$.
For any $k\in\{ 1,\dots ,m\}$, define
$$
G_k(z):=(-1)^k\rho (k-1)!\frac{L_k\left( z^{-\lambda _2}\right)}{z^k} ,
$$
where $z^{-\lambda _2}$ denotes any determination of the logarithm.
For sufficiently large $|z|$, this formula defines an analytic function on a slit plane.
Reporting in formula~(\ref{systemeFormel}) shows that the functions $G_k$ satisfy the simple differential
system
\begin{equation*}
%\label{systemeFormelSimple}
\left\{
\begin{array}{l}
\displaystyle
\forall k\in\{1,\dots ,m-2\},~
G'_k=G_{k+1},
\\
\displaystyle
G'_{m-1}=G_1^m.
\end{array}
\right.
\end{equation*}
In particular, $G_1$ is solution of the differential equation $y^{(m-1)}=y^m$.
We sum up these results in the following statement.

\begin{Th}
\label{thSerieLaplace}
Let $W_1$ be the complex-valued limit distribution for the multitype branching process of $m$-ary search trees
as defined in Section~\ref{fixedW}. Then:

\noindent
{\it (i)}
the Laplace series $L_1(z)=\g E\left( e^{zW_1}\right)$ has a positive radius of convergence;

\noindent
{\it (ii)}
for any determination of the logarithm, the function
$$
z\mapsto -\frac\rho zL_1\left( z^{-\lambda _2}\right) ,
$$
is a solution of the differential equation
\begin{equation}
\label{ym-1=ym}
y^{(m-1)}=y^m.
\end{equation}
\end{Th}

\begin{Rem}
As can be straightforwardly checked, the function $y_\kappa (z):=\frac{\kappa}{1-z}$ is a solution of
Eq.~(\ref{ym-1=ym}) when the complex number $\kappa$ satisfies $\kappa ^{m-1}=(m-1)!$.
Nonetheless, $G_1$ is not a function of this form.

Indeed, since $L_1(w)=1+w+o(w)$ in a neighbourhood of the origin, $G_1$ admits the expansion
$$
G_1(z)=-\frac\rho z-\frac{\rho}{z^{1+\lambda _2}}+o\left( \frac{1}{z^{1+\lambda _2}}\right),
$$
while $y_\kappa$ satisfies
$$
\frac{\kappa}{1-z}=-\frac\kappa z-\frac{\kappa}{z^2}+o\left( \frac 1{z^2}\right).
$$
One concludes by unicity of (complex) power expansions, because $\lambda _2\neq 1$.
\end{Rem}

\vskip 20pt\noindent
{\bf Acknowledgements.}
The authors owe much to Philippe Flajolet, especially some crucial arguments and many enthusiastic discussions. 
 %B+N : petit ajout de especially

%%%%%%%%%%%%%%%%%%%%%%%%%%%%%%%%%%%%%%%%%%%%%
\bibliographystyle{plain}
\bibliography{CLP}
%%%%%%%%%%%%%%%%%%%%%%%%%%%%%%%%%%%%%%%%%%%%

\end{document}